\documentclass[]{amsart}
\usepackage{amsfonts,amssymb,epsfig,latexsym}
\usepackage[latin1]{inputenc}
\usepackage[all]{xy}
\usepackage{amsmath}
\usepackage{delarray}
\usepackage{amssymb}
\usepackage{amscd}
\usepackage[mathscr]{eucal}
\usepackage{diagrams}

\usepackage{amsxtra}

\theoremstyle{plain}
\newtheorem{thm}{Theorem}[section]
\newtheorem{prop}{Proposition}[section]
\newtheorem{propdef}[prop]{Proposition--Definition}
\newtheorem{thdef}[thm]{Theorem--Definition}

\theoremstyle{definition}
\newtheorem{defn}{Definition}[section]
\newtheorem{ex}{Example}[section]

\newtheorem{conv}{Convention}[section]

\newtheorem{rem}{Remark}[section]

\numberwithin{equation}{section}

\def\textsl{\mathscr}
\def\textit{\mathcal}
\def\texttt{\mathbb}
\def\emphb{\boxed}

\def\textbfm{\mathbf}
\def\textbfs{\boldsymbol}

\newcommand{\lab}{labelstyle}

\newcommand{\sst}{\scriptstyle}

\newcommand{\bdi}{\begin{diagram}}
\newcommand{\edi}{\end{diagram}}

\newarrow{TO}----{curlyvee}

\newarrow{TeXto}----{->}
\newarrow{ClopTo}{sqhook}---{vee}
\newarrow{OpTo}{hooka}---{vee}
\newarrow{iTo}{triangle}---{curlyvee}

\newarrow{ITo}{blacktriangle}---{curlyvee}
\newarrow{sTo}----{triangle}
\newarrow{STo}----{blacktriangle}

\newarrow{Line}-----%

\newarrow{iLine}{triangle}----
\newarrow{ILine}{blacktriangle}----
\def\rH{\Rightarrow}
\def\lH{\Leftarrow}
\def\uH{\Uparrow}
\def\dH{\Downarrow}
\newarrowhead{=>}{\rH}{\lH}{\dH}{\uH}
\newarrow{DoubleLine}=====%
\newarrow{DoubleTo}===={=>}%
\newarrow{DoubleTO}===={curlyvee}%
\newarrow{Equal}=====%
\newarrow{DoublesTo}===={triangle}
\newarrow{DoubleSTo}===={blacktriangle}
\newarrow{DoubleiTo}{triangle}==={=>}
\newarrow{DoubleITo}{blacktriangle}==={=>}
\newarrow{DotTo}....{->}
\newarrow{DotsTo}....{triangle}
\newarrow{DotSTo}....{blacktriangle}
\newarrow{DotiTo}{triangle}...{->}
\newarrow{DotITo}{blacktriangle}...{->}
\newarrow{DashTO}{}{dash}{}{dash}{curlyvee}
\newarrow{DashTo}{}{dash}{}{dash}{vee}

\newarrow{DashsTo}{dash}{dash}{dash}{dash}{triangle}
\newarrow{DashSTo}{dash}{dash}{dash}{dash}{blacktriangle}
\newarrow{DashiTo}{triangle}{dash}{}{dash}{>}
\newarrow{DashITo}{blacktriangle}{dash}{dash}{dash}{>}
\newarrow{DashLine}{dash}{dash}{dash}{dash}{dash}
\newarrow{DashiLine}{triangle}{dash}{dash}{dash}{dash}
\newarrow{DashILine}{blacktriangle}{dash}{dash}{dash}{dash}
\newarrow{DotLine}.....

\def\rsum{\raisebox{2pt}{$\:\scriptscriptstyle+\:$}}
\def\lsum{\raisebox{2pt}{$\:\scriptscriptstyle+\:$}}
\def\dsum{\vphantom b\scriptscriptstyle+}
\def\usum{\vphantom b\scriptscriptstyle+}

\newarrowfiller{+}\rsum\lsum\dsum\usum
\newarrow{SumSTo}{+}{+}{+}{+}{blacktriangle}
\newarrow{SumSTo}{+}{+}{+}{+}{blacktriangle}
\newarrow{SumLine}{+}{+}{+}{+}{+}
\newarrow{SumTo}{+}{+}{+}{+}{->}
\newarrow{SumiTo}{triangle}{+}{+}{+}{->}

\def\rrnd{\raisebox{1.5pt}{$\:\scriptscriptstyle{\circ}\:$}}
\def\lrnd{\raisebox{1.5pt}{$\:\scriptscriptstyle{\circ}\:$}}
\def\drnd{\vphantom b\scriptscriptstyle{\circ}}
\def\urnd{\vphantom b\scriptscriptstyle{\circ}}

\newarrowfiller{rnd}{\rrnd}{\lrnd}{\drnd}{\urnd}
\newarrow{RndSTo}{rnd}{rnd}{rnd}{rnd}{blacktriangle}
\newarrow{RndLine}{rnd}{rnd}{rnd}{rnd}{rnd}

\def\rbul{\raisebox{1.5pt}{$\:\scriptscriptstyle{\bullet}\:$}}
\def\lbul{\raisebox{1.5pt}{$\:\scriptscriptstyle{\bullet}\:$}}
\def\dbul{\vphantom b\scriptscriptstyle{\bullet}}
\def\ubul{\vphantom b\scriptscriptstyle{\bullet}}

\newarrowfiller{bul}{\rbul}{\lbul}{\dbul}{\ubul}
\newarrow{BulSTo}{bul}{bul}{bul}{bul}{blacktriangle}
\newarrow{BulLine}{bul}{bul}{bul}{bul}{bul}

\def\imm{\text{\scriptsize{\textcircled{i}}}}
\def\rimm{\imm}
\def\limm{\imm}
\def\uimm{\imm}
\def\dimm{\imm}
\newarrowmiddle{imm}{\rimm}{\limm}{\uimm}{\dimm}
\newarrow{imTo}--{imm}-{->}
\newarrow{plfTo}{LaTeX}-{imm}-{->}
\newarrow{fidTo}{curlyvee}-{imm}-{->}
\newarrow{istTo}{boldlittlevee}-{imm}-{->}

\newarrow{qimTo}{o}---{->}
\def\etl{\text{\scriptsize{\textcircled{e}}}}
\def\retl{\etl}
\def\letl{\etl}
\def\uetl{\etl}
\def\detl{\etl}
\newarrowmiddle{etl}{\retl}{\letl}{\uetl}{\detl}
\newarrow{etaTo}--{etl}-{->}
\newarrow{etasTo}--{etl}-{triangle}
\newarrow{ploTo}{triangle}-{etl}-{->}
\newarrow{pofTo}{littleblack}-{etl}-{->}

\def\sub{\text{\scriptsize{\textcircled{s}}}}
\def\rsub{\sub}
\def\lsub{\sub}
\def\usub{\sub}
\def\dsub{\sub}
\newarrowmiddle{sub}{\rsub}{\lsub}{\usub}{\dsub}
\newarrow{subTo}--{sub}-{->}

\def\rph{\phantom{-}}
\def\lph{\phantom{-}}
\def\uph{\phantom{-}}
\def\dph{\phantom{-}}

\newarrowfiller{PH}{\rph}{\lph}{\uph}{\dph}
\newarrow{PHL}{PH}{PH}{PH}{PH}{PH}

\def\trl{\pitchfork}
\def\tre{\top}

\def\sq{\square\ }

\def\pb{\text{\raisebox{1.6pt}{{\fboxsep0.7pt\framebox{\tiny{\sffamily pb}}}}}}

\def\pbk{\SEpbk}

\def\ov{\overset}
\def\un{\underset}

\def\h{\quad}
\def\H{\qquad}
\def\HH{\H\H}
\def\HHH{\HH\H}
\def\HHHH{\HHH\H}
\def\HHHHH{\HHHH\H}
\def\v{\smallskip}
\def\V{\bigskip}

\usepackage{graphics}
\usepackage{color}

\definecolor{grreen}{cmyk}{0.8,0,1,0.5}

\def\emphb{\textcolor{blue}}
\def\emphr{\textcolor{red}}
\def\emphv{\textcolor{green}}

\definecolor{orange}{cmyk}{0,0.6,0.8,0.5}
\definecolor{violet}{rgb}{0.5,0,0.5}
\def\emphvi{\textcolor{violet}}
\def\emphor{\textcolor{orange}}

\begin{document}

\title[fractions]{Morphisms between spaces of leaves \\viewed as fractions}

\author{Jean PRADINES}
\address{26, rue Alexandre Ducos, F31500 Toulouse, France}
\email{jpradines@wanadoo.fr} 


\thanks{The author is indebted to Paul Taylor, whose package \emph{diagrams} is used throughout this text for drawing a variety of diagrams and arrows.}

\date{}

\begin{abstract}

Apr\`{e}s avoir transf\'{e}r\'{e}
 au cadre diff\'{e}rentiable la notion alg\'{e}brique d'\'{e}quivalence de groupo\"{i}des, nous montrons que les morphismes de la cat\'{e}gorie de fractions correspondante sont repr\'{e}sent\'{e}s par une unique fraction irr\'{e}\-duc\-tible (calcul de fractions simplifi\'{e}) que nous identifions aux morphismes de Connes--Skandalis--Haefliger entre espace de feuilles. Dans cette cat\'{e}gorie de fractions, le groupe fondamental de l'espace d'orbites au sens de Haefliger--van Est s'interpr\`{e}te comme r\'{e}flecteur sur la sous-cat\'{e}gorie pleine des groupes discrets.

The algebraic notion of equivalence between two groupoids is translated in the differentiable framework. Then we show that a morphism of the category of fractions in which these smooth equivalences are formally inverted may be represented by a unique irreducible fraction (simplified calculus of fractions, as in the elementary case of integers) which moreover may be identified with a generalized morphism in the sense of Connes--Skandalis--Haefliger between spaces of leaves. In this category of fractions, the fundamental group of the space of orbits in the sense of Haefliger--van Est is interpreted as a reflector onto the full sub-category of discrete groups.

\end{abstract}
\maketitle 

\setcounter{tocdepth}{2}
\tableofcontents

\begin{center}
0. \textsc{Introduction}
\end{center}
\V

The basic references for the present text are the papers by W.T. van Est \cite{vE} and A. Haefliger \cite{H}, in which various approaches to the transverse structure of foliations are described and certain concepts of transverse morphisms are introduced. The second approach is more general in that it considers topological groupoids which may be unequivalent to pseudogroups.

A very careful scrutiny of these papers would show that (when restricted to the common case of pseudogroups) the notions of morphisms considered by these authors are \emph{not} equivalent in general, though they are in the special cases of submersive morphisms and equivalences. More recently the ``generalized morphisms'' of A. Haefliger, attributed to G. Skandalis, have been used extensively, under the name of ``K-oriented morphisms'' by Skandalis, Hilsum \cite{HS} 
and the school of A. Connes \cite{AC}.

Here we start with van Est's \emph{geometrical} approach of pseudo-groups viewed as ``generalized atlases'', but we extend this (very illuminating) geometrical language to the ``non-\'{e}tale'' case, considering general groupoids as ``non-\'{e}tale atlases''. In this framework, a ``non-\'{e}tale change of base'' is an induction (or pullback) along a surmersion, which is a special case of equivalence (but will turn out to generate the most general concept).

It is then natural, from an \emph{algebraic} point of view, to define morphisms by formally inverting these surmersive equivalences, which is always possible in an abstract non-sense way \cite{GZ}. However the conditions for the classical calculus of fractions \cite{GZ} are \emph{not} fulfilled, but it turns out that we are able to unfold a ``simplified'' calculus of fractions in the sense of that our fractions admit unique irreducible representatives, as in the elementary case of integers. Now we have the remarkable fact that the irreducible fractions can be identified (in a non-obvious way) with the Skandalis-Haefliger morphisms.

The consideration of possibly non-reduced representatives gives a significant increase in flexibility. For instance the composition of morphisms becomes a routine diagram chasing (note that in the locally trivial topological case considered in \cite{H}, this composition is defined but in very special cases).

The irreducible fractions may also be viewed as special cases of J. B\'{e}nabou's distributors or profunctors (a more symmetrical notion). However the intersection of the two theories reduces to a rather trivial part of each one, and we let it to the informed reader \cite{BEN, J}.

 As an illustration we give a very simple characterization of the fundamental group of a foliation ( in the sense of van Est and Haefliger) by means of a reflection of our category of fractions into the full subcategory of discrete groups.

The present paper gives, essentially, ideas and results without detailed proofs. Our general policy throughout will be first to describe algebraic set-theoretic constructions by means of suitable diagrams in which we stress the injections and surjections, and secondly to replace injections by (regular) embeddings, and surjections by surmersions (i.e. surjective submersions). Then the proofs work by diagram chasing, using the formal properties of embeddings and surmersions listed in \cite{P75} under the name of ``diptych'' and the formal properties of commutative squares stated in the basic proposition A2 of \cite{P85}.

In the following, a \emph{pseudogroup} of transformations (always assumed to be complete or completed) will always be identified with the \emph{groupoid} of its germs, provided with the (\'{e}tale) sheaf topology.

\section{The language of (generalized) atlases.}
\label{sec:language of atlases}

Let us first consider a (smooth) manifold $Q$ and a (classical) atlas of $Q$, i.e. a collection of charts $p_i:V_i\rightarrow U_i$ (open sets in some $\texttt{R}^n$), or equivalently of cocharts $q_i=p_i^{-1}$. It is equivalent to consider the \'{e}tale surjective map $q:U\rightarrow Q$ where $U$ is the (trivial) manifold coproduct (or disjoint sum) of the $U_i$'s. The fibred product $R=U\un{Q}{\times}U$, with its projections $\alpha=\text{pr}_2$, $\beta=\text{pr}_1$, may be viewed either as the graph of the equivalence relation in $U$ defined by $q$ or as the pseudogroup of changes of charts, which is a (very special kind of) groupoid with base $U$. Conversely the data of $R$ with its manifold and groupoid structures determine $Q$ and $q$ up to isomorphisms.
 
 In that context, a refinement of the given atlas is viewed as an \'{e}tale surjective map $u: U'\rightarrow U$ and then the corresponding graph $R'$ is obtained by pulling back along $u$. Two atlases are equivalent if they admit a common refinement.

This situation admits a twofold generalization.

First following van Est a pseudogroup may be viewed as a generalized (\'{e}tale) atlas of its space of orbits (which is no longer a manifold in general). This applies to any regular foliation, using a totally transverse manifold $T$ 	and the corresponding holonomy pseudogroup, whose space of orbits is the space of leaves. Various choices of $T$ lead to equivalent atlases in a generalized sense explained below.

Second replacing $q$ by a (possibly non-\'{e}tale) surmersion $q:B\rightarrow Q$, we can view the graph $R=B\un{Q}{\times}B$ ( with its manifold and groupoid structures) as a ``non-\'{e}tale atlas'' of $Q$ with base $B$. A non-\'{e}tale refinement is then a surmersion $B'\rightarrow B$ and the new ``atlas'' $R'$ is again obtained by pulling back. If moreover  $q$ is ``retro-connected'' (i.e. its fibres are connected), the manifold $Q$ is the space of leaves of the simple foliation of $B$ defined by $q$.
 
 A further generalization is required for a non-simple (regular) foliation, the previous construction being valid only locally. The local pieces can be glued together into the holonomy groupoid introduced by Ehresmann in \cite{OC}/54/ (see also \cite{P66, P84}) and renamed graph of the foliation by Winkelnkemper \cite{W} and A. Connes \cite{AC}. Though this groupoid has special properties which we emphasized in \cite{P84}, we do not use them in the sequel.
So we are led to the following common generalization. This generalization makes use of the general notion of smooth (or differentiable) groupoid introduced by Ehresmann (\cite{OC}/50/) which we recall first.

\section{Smooth groupoids and orbital atlases.}
\label{orbat}
In the sequel $\textbfs{\textsl{D}}$ will denote the category of (morphisms between) smooth manifolds. We consider the following subcategories:
\begin{itemize}
	\item $\textbfs{\textsl{D}}^\ast = \text{diffeomorphisms}$\,;
	\item $\textbfs{\textsl{D}}_e = \text{\'{e}tale maps ( or local diffeomorphisms)}$\,;
	\item $\textbfs{\textsl{D}}_i = \text{(regular) embeddings, denoted by } \riTo$\,;
	\item $\textbfs{\textsl{D}}_s = \text{surmersions, denoted by }\rsTo$\,;
	\item $\textbfs{\textsl{D}}_{ei}=\textbfs{\textsl{D}}_e\cap\textbfs{\textsl{D}}_i$\,;\h$\textbfs{\textsl{D}}_{es}=\textbfs{\textsl{D}}_e\cap\textbfs{\textsl{D}}_s$\,.
\end{itemize}

The subclass $\textbfs{\textsl{D}}_i\textbfs{\textsl{D}}_s$ (which is not a subcategory !) is denoted by $\textbfs{\textsl{D}}_r$ (regular morphisms).

Let \textbf{P}:
\begin{align}
	\bdi[h=2.5em,w=3em,tight,\lab=\sst]
A'&\rTo^{f'}&B'\\
\dTo^u&\textbf{P}&\dTo_v\\
A&\rTo^f&B
\edi\notag
\end{align}
be a commutative square of $\textbfs{\textsl{D}}$, and denote by $R=A\un{B}{\times}B'$ the (set-theoretic) fibred product of $A$ and $B'$ over $B$. Then $\textbf{P}$ is called:
\begin{itemize}
	\item \emph{i-faithful} if $(u,f'):A'\rightarrow A\times B'$ lies in $\textbfs{\textsl{D}}_i$\,;
	\item \emph{universal} (resp. \emph{s-full}) if $R$ is a submanifold of $A\times B'$ and moreover the canonical map $A'\rightarrow R$ lies in $\textbfs{\textsl{D}}^\ast$ (resp. $\textbfs{\textsl{D}}_s$).
\end{itemize}

Note that universal implies $\textbfs{\textsl{D}}$-cartesian (i.e. pullback square in $\textbfs{\textsl{D}})$ but \emph{the converse is false}. Note also that the transversality \cite{VAR} of $f$ and $v$ \emph{implies} that the pullback is universal but \emph{the converse is false}: we shall say that $f$ and $v$ are \emph{weakly transversal} when they can be completed into a universal square in the sense just defined.

The basic properties of such squares are stated (with a different terminology) in Proposition A2 of \cite{P85}, which we complete by the following: 
\begin{center}
\emph{If $f$ is a surmersion, and $\textbfm{P}$ and $\textbfm{QP}$ are universal, then $\textbfm{Q}$ is universal}.
\end{center}

Now we remind that a (small) \emph{groupoid} is a (small) category with all arrows invertible. Usually a groupoid will be loosely denoted by its set of arrows $G$.
The base $B=G_0$ is the set of \emph{objects}, identified by the \emph{unit map} $\omega_G: B\rightarrow G$ with the set of \emph{units} $\omega_G(B)\subset G$. The \emph{source} and \emph{target} maps are denoted by: $\alpha_G$, $\beta_G:G\rightarrow B$. The map:
	\[\tau_G=(\beta_G,\alpha_G):G\rightarrow B\times B
\]
will be called the \emph{transitor} (anchor map in \cite{MK}). The image of $\tau_G$ is the graph of an equivalence relation in $B$ whose classes are the \emph{orbits} of $G$ in $B$. The inverse images of the orbits are the \emph{transitive components} of $G$.
The map
	\[\delta_G: \Delta G\rightarrow G,\h (x,y)\mapsto xy^{-1}
\]
(where $ \Delta G\subset G\times G$ is the set of pairs of arrows with the same source)
 may be called the \emph{divisor}.

The morphisms $f:G'\rightarrow G$ between groupoids are just the \emph{functors} and are the arrows of a category $\textbfs{\textsl{G}}$. The restriction $f_0:B'\rightarrow B$ of $f$ to the bases of $G'$, $G$ may be called the \emph{objector} of $f$ ; when $f_0$ is the identity of $B$, $f$ is said to be \emph{uniferous}. The subcategory of uniferous functors will be denoted by $\textbfs{\textsl{G}}_0$, and $\textbfs{\textsl{G}}_B$ when the base $B$ is fixed.

We say that the \emph{groupoid} $G$ is \emph{smooth} (or \emph{differentiable}) (\cite{OC}/50/, \cite{MK}) when $G$ and $B$ are provided with manifold structures such that $\omega_G\in\textbfs{\textsl{D}}$, $\alpha_G\in \textbfs{\textsl{D}}_s$ (which implies that $\Delta G$ is a submanifold of $G\times G$), and $\delta_G\in \textbfs{\textsl{D}}$. This implies easily $ \omega_G \in \textbfs{\textsl{D}}_i$, $\beta\in \textbfs{\textsl{D}}_s$, $\delta_G\in \textbfs{\textsl{D}}_s$.

A \emph{functor} $f:G\rightarrow G'$ is \emph{smooth} if the underlying map is smooth: if moreover it lies in $\textbfs{\textsl{D}}_i$ (resp. $\textbfs{\textsl{D}}_s$), we say that $f$ is an \emph{i-functor} (resp. \emph{s-functor}): note that this implies $f_0$ is also in $\textbfs{\textsl{D}}_i$ (resp. $\textbfs{\textsl{D}}_s$). The category of smooth functors between smooth groupoids is denoted by $\textbfs{\textsl{GD}}$.

A smooth functor is \emph{split} when it admits a section in $\textbfs{\textsl{GD}}$.

To any smooth functor $f: G\rightarrow G'$ there are associated two commutative squares:

\begin{align}\label{AT}
	\begin{diagram}[h=3em,w=4em,tight,inline,\lab=\sst]
G'&\rTo^{f}&G\\
\dsTo^{\alpha_{G'}}&\textbf{A}(f)&\dsTo_{\alpha_G}\\
B'&\rTo^{f_0}&B
\end{diagram}
\HH
\begin{diagram}[h=3em,w=4em,tight,inline,\lab=\sst]
G'&\rTo^{f}&G\\
\dTo^{\tau_{G'}}&\textbf{T}(f)&\dTo_{\tau_G}\\
B'\times B'&\rTo^{f_0\times f_0}&B\times B \tag{A, T}
\end{diagram}
\end{align}
the first one in $\textbfs{\textsl{D}}$, the second one in $\textbfs{\textsl{GD}}$.

A smooth functor is called \emph{i-faithful} (\emph{s-full}, an \emph{inductor}) when the square $\textbf{T}(f)$ is $i$-faithful ($s$-faithful, universal). From Proposition A2 of \cite{P85} we get:
\begin{prop}
Let $h=gf$ be the composite of two smooth functors.
\begin{enumerate}
	\item[(i)] If $f$ and $g$ are i-faithful (resp. inductors), so is $h$;
	\item[(ii)] If $h$ is i-faithful, so is $f$.
	\item[(iii)] Assume $f$ is an s-functor and an inductor (briefly an s-inductor); then if $h$ is i-faithful (an inductor), so is $g$;
	\item[(iv)] Assume $g$ is an inductor; then $f$ is s-full (an inductor) iff $h$ is.
\end{enumerate}
\end{prop}

Now the considerations of section 1 lead us to set:
\begin{defn} An \emph{orbital atlas} on a set $Q$ is a pair $(G,\,q)$ where $G$ is a smooth groupoid with base $B$ and $q:B\rightarrow Q$ is a surjection whose fibres are the orbits of $G$ in $B$.
\end{defn}

$Q$ will be provided with the finest topology making $q$ continuous. Then $q$ is \emph{open}.

A basic example is the holonomy groupoid of a foliation, viewed as an orbital atlas of the space of leaves. 

Note that transitive smooth groupoids (particularly Lie groups) define various unequivalent orbital atlases for a singleton.

\section{Surmersive equivalences and extensors.}
\label{sec:SEquivalencesExtensors}

If $u:B'\rsTo B$ lies in $\textbfs{\textsl{D}}_s$, the fibred product $G'=u^\ast(G)$ of the arrows $\tau_G$ and $u\times u$ has a canonical structure of groupoid called the pullback of $G$ along $u$, for which $f:G'\rightarrow G$ is an $s$-inductor. Any smooth functor $g:H\rightarrow G$ with its objector $g_0=u$ admits a unique factorization $g=fh$. 

\begin{defn} An \emph{s-inductor} will be called also an \emph{s-equivalence}; an \emph{s-extensor} is an $s$-full $s$-functor.
\end{defn}

The following statements are proved in \cite{P86}:
\begin{thm}
\begin{itemize}
	\item[(i)]An s-equivalence induces an equivalence between the categories $(G\downarrow \textbfs{\textsl{GD}}_B)$ and $(G'\downarrow \textbfs{\textsl{GD}}_{B'})$ of groupoids under $G$ and  $G'$ \cite{ML}.
	\item[(ii)] Let $f:H\rightarrow G$ be a smooth functor and $N=f^{-1}(B)$ its set-theoretic kernel: then the following statements are equivalent:
	
\begin{itemize}
	\item[a)] $f$ is an extensor;
	\item[b)] $N$ is a regular smooth groupoid embedded in $H$ and the square:
	
\begin{center}
\bdi[h=2.5em,w=3
em,tight,inline,\lab=\sst]
N\pbk &\rTo&B\\
\diTo&&\diTo\\
H&\rTo^f&\NWpbk G\\
\edi
\end{center}
is a pushout in $\textbfs{\textsl{GD}}$;
\item[c)] $N$ is a regular smooth groupoid embedded in $H$, $f$ is an s-functor, and the relation $f(x)=f(y)$ is equivalent to $x\in NyN$ (two-sided coset).
\end{itemize}
\end{itemize}
\end{thm}

Keeping the above notations, if $(G,\,q)$ is an orbital atlas of $Q$, then $(G',\,q')$, where $q'=uq$ (with $u$ an s-equivalence) is again an orbital atlas of $Q$ called a \emph{refinement} of $(G,\,q)$. Two atlases of $Q$ are said to be \emph{equivalent} if they admit a common refinement.

It is convenient to think an equivalence class of orbital atlases on $Q$ as defining a (generalized!) ``structure'' on the the set $Q$, called \emph{orbital structure}. But one should notice carefully that the morphisms we shall introduce will be defined only at the atlas level and not between such ``structures'', which do not play the role of objects of some category.

Two smooth groupoids $G_i$, $(i=1,\,2)$ are called (smoothly) \emph{equivalent} if there exists a pair of $s$-equivalences $f_i:G\rightarrow G_i$ : this is indeed an equivalence relation.

\begin{prop}
Let $h=gf$ be the composite of two smooth functors. Then:
\begin{itemize}
	\item[(i)] if $f$ and $g$ are s-extensors, so is $h$;
	\item[(ii)] assume $g$ is an s-equivalence and $f_0\in \textbfs{\textit{D}}_s$\,: then if $h$ is an s-extensor or an s-equivalence, so is $f$;
	\item[(iii)] assume $f$ is an s-extensor: then if $h$ is an s-extensor or an s-equivalence, so is $g$.
\end{itemize}
\end{prop}

\section{Some important special smooth groupoids.}
\label{sec:SpecialGroupoids}
Let $G$ be a smooth groupoid with base $B$. We consider various special cases.
\begin{itemize}
	\item[(i)] $G$ is (topologically!) \emph{discrete}; we can identify $\textbfs{\textsl{G}}$ with the full subcategory of topologically discrete smooth groupoids in $\textbfs{\textsl{GD}}$;
	\item[(ii)] $B$ is a \emph{singleton}: $G$ is (identified with) a \emph{Lie group}.
	\item[(iii)] $\alpha_G=\beta_G$\,: $G$ is called a smooth \emph{plurigroup}: the full subcategory of smooth (pluri)\-groups will be denoted by $\textbfs{g}\textbfs{\textsl{D}}$\h(\,$\overline{\textbfs{g}}\,\textbfs{\textsl{D}}$\,);
	\item[(iv)] $\omega_G\in \textbfs{\textsl{D}}^\ast$\,: $G$ is \emph{null}; we may identify $\textbfs{\textsl{D}}$ with the full subcategory of null smooth groupoids in $\textbfs{\textsl{GD}}$\,;
	\item[(v)] $\tau_G\in \textbfs{\textsl{D}}^\ast$\,: $G$ is \emph{banal};
	\item[(vi)] $\tau_G\in \textbfs{\textsl{D}}_i$\,: $G$ is \emph{principal} (or \emph{Godement}); by Godement's Theorem, $G$ is identified with the graph of a regular equivalence relation in $B$;
	\item[(vii)] $\tau_G\in \textbfs{\textsl{D}}_s$\,: $G$ is \emph{s-transitive}; the fibres of $\alpha_G$ are principal bundles with base ~$B$, and $G$ may be identified with their structural groupoid \cite{OC}/50/;
	\item[(viii)] $\tau_G\in \textbfs{\textsl{D}}_{es}$\,: $G$ is \emph{es-transitive} or a \emph{Galois} groupoid (structural groupoid of a Galois or normal covering);
	\item[(ix)]  $\tau_G\in \textbfs{\textsl{D}}_r$\,: $G$ is \emph{regular};
	\item[(x)] $\tau_G$ is a weak embedding: $G$ is a \emph{Barre} groupoid (its space of orbits is a Barre Q-manifold) \cite{BAR};
	\item[(xi)] $\tau_G$ is a faithful immersion: $G$ is a \emph{graphoid} in the sense of \cite{P84}.
\end{itemize}

\begin{prop} $G$ is principal (s-transitive, a Galois groupoid, a graphoid) iff it is equivalent to a null groupoid (Lie group, discrete group, pseudo\-group).
\end{prop}

The holonomy groupoid of a regular foliation is equivalent to any of its transverse holonomy pseudo\-groups.

\begin{defn}
A smooth functor is called \emph{principal} if its source is principal.
\end{defn}

\begin{prop}
Assume the smooth functor $f:H\rightarrow G$ is i-faithful (resp. s-full, resp. an s-extensor)). Then if $G$ is principal (resp. Lie, resp. regular), so is $H$.
\end{prop}

\section{Smooth equivalences.}
\label{sec:SmoothEquivalences}

Following our general policy, we give a smooth version of the algebraic notions of essential (or generic) surjectivity and equivalences between groupoids (more general than the surjective equivalences).

Let be given a smooth groupoid $G$ with base $B$ and a map $b:B'\rightarrow B$ (in $\textbfs{\textsl{D}}$). Let $W$ be the fibred product (in $\textbfs{\textsl{D}}$) of $\alpha_G$ and $b$ and consider the following diagram in $\textbfs{\textsl{D}}$:
\begin{center}
\bdi[h=1.4em,w=1.8em,tight,inline,\lab=\sst]
&{}&&&\rLine^v[abut]&&&{}&\\
W\pbk\ruLine(1,1)[abut]&&\rTo_u&&G&&\rsTo_{\beta_G}&&B\rdTo(1,1)[abut]\\
&&&&&&&&\\
\dsTo^a&&&&\dsTo_{\alpha_G}\\
&&&&&&&&\\
B'&&\rTo^b&&B&&&&\\
\edi
\end{center}

\begin{defn}
\label{essur}
We say $b$ is \emph{transversal} to $G$ when $v$ lies in $\textbfs{\textsl{D}}_s$, and a functor $f:G'\rightarrow G$ is \emph{essentially surmersive} when $f_0$ is transversal to $G$.
\end{defn}

\begin{propdef}\label{ind}
If $b$ is transversal to $G$, then the fibred product of $b\times b$ and $\tau_G$ does exist in $\textbfs{\textsl{GD}}$ and the pullback we get is universal in $\textbfs{\textsl{D}}$. We say that $G'$ is the (smooth) groupoid \emph{induced} by $G$ along $b$ (or the \emph{pullback} of $G$ along ~$b$).

A smooth functor $f:H\rightarrow G$ with $f_0=b$ is called a (\emph{smooth}) \emph{equivalence} if it is essentially surmersive and if the canonical factorization $H\rightarrow G'$ is an isomorphism.
\end{propdef}

\begin{prop}
\begin{itemize}
	
	\item[(i)] The equivalences and the essentially surmersive functors make up subcategories of $\textbfs{\textsl{GD}}$.
	\item[(ii)] If $g$ is an equivalence and $gf$ is essentially surmersive (resp. is an equivalence), then $f$ is essentially surmersive (resp. is an equivalence).
	\item[(iii)]If $f_0$ lies in $\textbfs{\textsl{D}}_s$ (resp. if $f$ is an s-extensor) and $gf$ is essentially surmersive (resp. is an equivalence), then $g$ is essentially surmersive (resp. is an equi\-valence and $f $ is an s-equivalence).
\end{itemize}
\end{prop}

\section{Holomorphisms.}
\label{sec:Holo}

If $\sq G$ denotes the smooth groupoid of \emph{commutative squares} of $G$ with the horizontal composition law, the two canonical projections $\varpi_1$, $\varpi_2$ onto $G$ are $s$\nobreakdash-equi\-va\-lences while the canonical injection $\iota_G$ is an $i$-equivalence (and a common section of $\varpi_1$, $\varpi_2$).

A (smooth) \emph{natural transformation} between two smooth functors $f_1,\,f_2:G\rightarrow H$ may be described either a as smooth functor $\texttt{I}\times G\rightarrow H$ (where $\texttt{I}$ is the banal groupoid $\{0,\,1\}\times\{0,\,1\}$) or a smooth functor $G\rightarrow \sq H$. As a consequence:
\begin{prop}
The following properties of a smooth functor are preserved by a smooth functorial isomorphism: $i$-faithful, $s$-full, essentially surmersive, equi\-valence, $s$-extensor.
\end{prop}

By the horizontal composition of natural transformations, the isomorphism between smooth functors is compatible with the composition of functors.

This gives rise to a new category (with the \emph{same objects} as $\textbfs{\textsl{GD}}$) denoted by $\boxed{[\textbfs{\textsl{G}}]\textbfs{\textsl{D}}}$, the arrows of which will be called \emph{holomorphisms}, and to a canonical full functor $f\mapsto [f]$ from $\textbfs{\textsl{GD}}$ to $[\textbfs{\textsl{G}}]\textbfs{\textsl{D}}$.

The holomorphisms between Lie groups  are just the \emph{conjugacy classes} of homomorphisms. So the notion of holomorphism extends the notion of outer automorphism (this suggests the alternative terminology of \emph{exomorphism}).

\section{Actors, exactors, subactors.}
\label{sec:ActorsExactorsSubactors}

After the diagram $\textbfm{T}(f)$, which measures the (lack of) faithfulness of $f$, we turn now to the diagram $\textbfm{A}(f)$, which measures its ``activity'' ( i.e. how far it is from describing an action law). (In the purely algebraic context several variants of the notions below have been used by various authors such as Ehresmann, Grothendieck, Higgins, R. Brown, van Est \emph{et alii}, under various names, notably (discrete) (op)fibrations, coverings, and others, which we cannot carry over to the smooth case.)

\begin{defn}
\label{act}
A smooth functor $f$ is called an \emph{actor} (\emph{inactor}, \emph{exactor}) when the square $\textbfm{A}(f)$ (section \ref{orbat}) is universal ($i$-faithful, $s$-full). More precisely we speak of $G$-actor, when the target $G$ is fixed.

\emph{The actors / exactors will be tagged by @ / $\widehat{@}$.}

There is an equivalence of categories between the category of (morphisms between) $G$-actors and the category of (equivariant morphisms between) smooth action laws of $G$ on manifolds over the base $B$ of $G$ (hence the terminology)\cite{MK}.
\end{defn}

\begin{rem}
\begin{itemize}
	\item[(i)] The image of an actor is a (possibly non-smooth) sub\-group\-oid of $G$.
	\item[(ii)] Any $s$-extensor is an $s$-exactor; any inactor is $i$-faithful.
	\item[(iii)] An exactor is essentially surmersive iff it is an $s$-exactor.
\end{itemize}
\end{rem}

\begin{prop}
A smooth functor which is an equivalence and an actor is an isomorphism (of smooth groupoids). If it is an exactor and an inductor, it is an $s$-equivalence.
\end{prop}

\begin{prop}
If $f:G'\rightarrow G$ is an $s$-exactor, $H$ a smooth groupoid, and $h:G\rightarrow H$ a (set-theoretic) map such that $hf: G'\rightarrow H$ is a smooth functor, then $h:G\rightarrow H$ is a smooth functor.
\end{prop}

\begin{prop}
Let $h=gf$ be the composite of two smooth functors.
\begin{itemize}
	\item[(i)] If $f$, $g$ are (ex)(in)actors, so is $h$.
	\item[(ii)] Assume $g$ is an actor. Then if $h$ is an (ex)actor, so is $f$.
	\item[(iii)] Assume $f$ is an $s$-exactor. Then if $h$ is an (ex)actor, so is $g$.
	\item[(iv)] If $h$ is an inactor, so is $f$.
\end{itemize}
\end{prop}

\begin{prop}\label{pb}
Let $g: G'\rightarrow G$ be an \emph{(ex)actor}, and $u:H\rightarrow G$ a smooth functor. Assume $g_0$ and $u_0$ to be \emph{weakly transversal}. Then:
\begin{itemize}
	\item[(i)] The \emph{fibred product} $H'=G'\un{G}{\times}H$ exists in $\textbfs{\textsl{GD}}$, the pullback square is \emph{universal} in $\textbfs{\textsl{D}}$, and $h:H'\rightarrow H$ is an \emph{(ex)actor}. The induced map \emph{$k:\text{Ker}\,h\rightarrow\text{Ker}\,g$} is an \emph{actor}.
	\item[(ii)] If moreover $g$ is an $s$-extensor (an $s$-equivalence), so is $h$.
	\item[(iii)]  If $u$ is an inactor ($i$-faithful, essentially surmersive, an inductor, an equivalence), so is $u':H'\rightarrow G'$. If moreover $g$ is an $s$-exactor, then if $u'$ is an (in)(ex)actor (essentially surmersive, an inductor, an equivalence), so is $u$.
\end{itemize}
\end{prop}

The situation is depicted by the following cubic diagram:
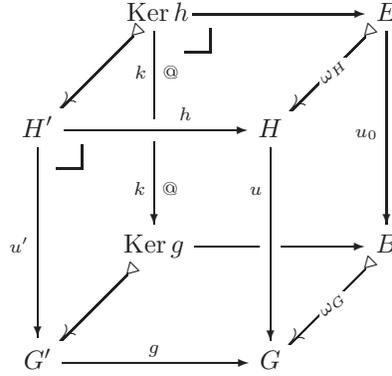
\begin{center}
\bdi[h=2.2em,w=2.2em,tight,inline,\lab=\sst]
&&{\text{Ker}\,h}\negthickspace\pbk&& \rTo&&E\\
&\ldiTo&\dLine^k_{@}&&&\ldiTo~{\omega_H}&\\
H'\pbk&&\HonV&\rTo^{h\H}&H&&\dTo^{u_0}\\
&&\dTo^k_{@}&&\dTo^{u}&&\\
\dTo^{u'}&&\text{Ker}\,g&\rLine&\VonH&\rTo&B\\
&\ldiTo&&&&\ldiTo~{\omega_G}&\\
G'&&\rTo^g&&G&&\\
\edi
\end{center}
in which the front, rear, top and bottom faces are pull back squares.

Observe that any exactor $f$ has a \emph{kernel} in $\textbfs{\textsl{GD}}$; $f$ will be an actor iff this kernel is null.

The more general case when this kernel is \emph{principal} is of importance too:

\begin{propdef}
Let $f:H\rightarrow G$ be an \emph{exactor}. The following are equivalent:
\begin{itemize}
	\item[(i)] \emph{$\text{Ker}f$} is principal \emph{(section \ref{sec:SpecialGroupoids}, special case (vi))};
	\item[(ii)] $f$ is $i$-faithful;
	\item[(iii)] $f=ae$ where $e$ is an $s$-equivalence
 and $a$ an actor.
 The decomposition \emph{(iii)} is essentially unique.
 \end{itemize}
 Then $f$ is called a \emph{subactor}.
\end{propdef}

\begin{rem}
It will be proved elsewhere that any $i$-faithful functor is the composite of an equivalence and an actor.
\end{rem}

The following two propositions generalize a lemma of van Est \cite{vE} (p. 245).

\begin{prop}
Assume $ae'=ea'$ where $a$, $a'$ are actors, $e'$an equivalence, and $e$ an $s$-equivalence. Then the square is a pullback (hence $e'$ is also an s-equivalence).
\end{prop}

This is displayed in the following diagram:
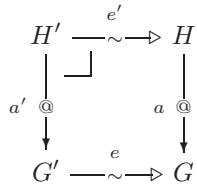
\begin{center}
\bdi[h=1.3em,w=1.3em,tight,inline,\lab=\sst]
H'\pbk&&\rsTo~\sim^{e'}&&H\\
&&&&\\
\dTo^{a'}~{@}&&&&\dTo^{a}~{@}\\
&&&&\\
G'&&\rsTo~\sim^{e}&&G\\
\edi
\end{center}

Now let $u:G'\rightarrow G$ be an $s$-equivalence. Then pulling back along $u$ determines a functor $u^\ast: (\text{Act}\downarrow G)\rightarrow(\text{Act}\downarrow G')$ from the category of $G$-actors to the category of $G'$-actors. Conversely we define the direct image of a $G'$-actor $a'$ by taking for $u_\ast(a')$ the first factor of the decomposition (iii) above.

\begin{thm}
The pair $(u^\ast,\,u_\ast)$ defines an adjoint equivalence \cite{ML} between \emph{$(\text{Act}\downarrow G)$} and \emph{($\text{Act}\downarrow G')$}.
\end{thm}

\section{Holograph of a functor}
\label{Hologr}

The following smooth construction is known in the algebraic context of profunctors \cite{BEN, J}. It turns out to be crucial for defining the (non-trivial) functor from the functors to the fractions.

Let $f:H\rightarrow G$ be an arrow of $\textbfs{\textsl{GD}}$.

Using Prop. \ref{pb} (ii) in order to pull $f$ back along $\varpi_2$ (notations of section \ref{sec:Holo}) we can construct the commutative diagram in $\textbfs{\textsl{GD}}$\,:
\begin{center}
\begin{align}
	\bdi[h=3em,w=3.5em,tight,inline]
	&& H^{\square}\negthickspace\negthickspace\pbk\h&&\rTo~{f^\square}&&\sq G&&\\
	&
\ldsTo~{\,q=q(f)\;\sim\,}(2,4)&&
\rdTo~{\,p=p(f)\;\widehat{@}\,}(6,4)
\rdTo~{\H}~{\,p'= p'(f)\;\widehat{@}\,}(2,4)&&
\ldsTo~{\,\varpi_2\H\sim\,}(2,4)&&
\rdsTo(2,4)~{\,\varpi_1\;\sim\,}&\\
	&&&
	&&&&&\\
	&&&&&&&&\\
	H&&\rTo~f&&G&&&&G\\
	\edi\tag{holo}
\end{align}
\end{center}
which associates to $f$ the pair of arrows $(p,\,q)$, with $p=p(f)=\varpi_1 f^\square$,\h $q=q(f)$. 
\begin{propdef}
For any smooth functor $f$, $p$ is an \emph{exactor} and $q$ a \emph{split $s$-equivalence}. We call $(p,\,q)$ the \emph{holograph} of $f$,\;$q=q(f)$ its \emph{denominator}, and $p=p(f)$ the \emph{expansion} \emph{(\emph{or} numerator)} of $f$. One observes that $p$ is \emph{isomorphic} to $p'=fq=\varpi_2f^\square$, by means of the natural transformation defined by $f^{\square}$.
\end{propdef}

The \emph{holograph of the identity} of $G$ is $(\varpi_1,\,\varpi_2)$.

\begin{prop}
A smooth functor $f$ is essentially surmersive (i-faithful, an equivalence) iff its expansion $p(f)$ is an s-exactor (a subactor, an s-equivalence).
\end{prop}

\begin{ex}
The \emph{holograph of the unit map} $\omega_G:B\riTo G$ is $(\delta_G,\, w_G)$ where $w_G$ is the canonical ``vertex'' projection onto the common source $w_G:\Delta_G\rsTo^{\sim\h} B$ (an $s$-equivalence); here the base $B$ of $G$ is regarded as a null groupoid.
\end{ex}

\section{Transversal and transverse subgroupoids}
\label{sec:TransversalSubgroupoids}

Let $K$ be a smooth groupoid with base $E$, and $M$, $N$ two uniferous embedded subgroupoids, $i$, $j$ their canonical injections, $S$ the (generally non-smooth) subgroupoid $M\cap N$.

Let $L$ be the fibred product of $\alpha_M$ and $\alpha_N$, which is a submanifold of $\Delta K$.

\begin{defn}
$M$ and $N$ are called \emph{transversal} in $K$ (denoted by $M\trl N$) if the restriction of $\delta_K$ to the submanifold $L$ is a \emph{surmersion} onto $K$. They are called \emph{transverse} ($M\tre N$) if it is a \emph{diffeomorphism} (then $S$ is null).
\end{defn}

When such is the case, it can be proved that $S$ is a \emph{smooth subgroupoid embedded} in $M$ and $N$; in particular, if $M$ or $N$ is \emph{principal}, so is $S$.

\begin{rem}
The data $M$, $N$ with $M\tre N$ determine on $K$ a structure of \emph{smooth double groupoid} (\cite{OC} /63/): $M$ and $N$ are the respective bases of the horizontal and vertical laws and the source map $K\rightarrow M$ of the horizontal law is an $s$-actor when $K$ and $M$
 are considered with the vertical law. The converse is true. We do not develop these facts that are not needed here.
\end{rem}

\begin{propdef}
Let be given an exactor $p:K\rightarrow G$, and assume \emph{$N=\text{Ker}\,p$}.
 Let $M$ be another uniferous subgroupoid embedded in $K$. Then one has $M\trl N$ (resp. $M\tre N$) iff $u=pi$ is an \emph{exactor }(resp. an \emph{actor});  when such an $M$ exists, we say $p$ is \emph{inessential}. \emph{(Note that for surjective homomorphisms of groups the notions of inessential and split coincide.)}
 \end{propdef}
 
 As a consequence, if $M$ is also the kernel of an exactor $q:K\rightarrow H$, then $u=pi$ is an (ex)actor iff $v=qj$ is. If such is the case we say the exactors $p$ and $q$ are \emph{cotransvers}(\emph{al}). The situation is pictured by the following butterfly diagram:
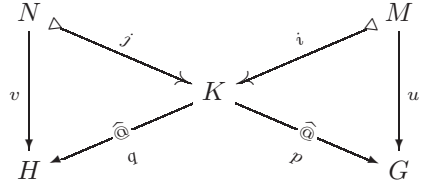
\begin{center}
\bdi[h=1.5em,w=3.5em,tight,inline,\lab=\sst]
N&&&&M\\
&\rdiTo^j(2,2)&&\ldiTo^i(2,2)&\\
\dTo^v&&K&&\dTo_{u}\\
&\ldTo_{q}~{\widehat{@}}&&\rdTo_{p}~{\widehat{@}}\\
H&&&&G\\
\edi
\end{center}
in which $p$, $q$ are \emph{exactors}, $N=\text{Ker}\,p$,\;$M=\text{Ker}\,q$; the condition of co\-transvers(al)ity is reflected by the property of $u$ and $v$ being both (ex)actors.

\section{Fractions and meromorphisms: the simplified calculus of fractions.}
\label{sec:Mero}
We consider now the category whose objects are pairs $(p,\,q)$ of \emph{exactors} with the same source $H\stackrel{q}{\leftarrow}K\stackrel{p}{\rightarrow }G$, and arrows $k:(p',\,q')\rightarrow(p,\,q)$ are smooth functors $k:K'\rightarrow K$ making the whole diagram commutative, i.e. $p'=pk,\h q'=qk$, as shown in the following diagram:
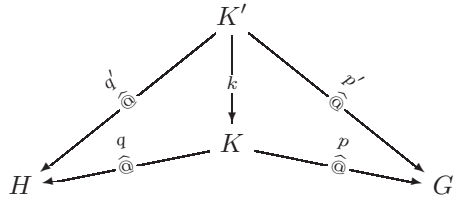
\begin{center}
\bdi[h=.8em,w=4em,tight,inline,\lab=\sst]
&&K'&&\\
&\ldTo^{q'}~{\widehat{@}}(2,8)&&\rdTo^{p'}~{\widehat{@}}(2,8)&\\
&&&&\\
&&\dTo~k&&\\
&&&&\\
&&&&\\
&&K&&\\
&\ldTo^q~{\widehat{@}}&&\rdTo^{p}~{\widehat{@}}&\\
H&&&&G\\
\edi
\end{center}

The \emph{isomorphy class} of the pair $(p,\,q)$ (this means $k\in\textbfs{\textsl{D}}^\ast$) will be denoted by $\boxed{p/q}$ and called a \emph{fraction} with source $H$ and target $G$.

Two pairs $H\stackrel{q_i}{\leftarrow}K_i\stackrel{p_i}{\rightarrow}G\h(i=1,\,2)$ are \emph{equivalent} if there exist two $s$-equivalences $K\rsTo~{\sim}^{k_i}K_i$ making the whole diagram commutative (this is indeed an equivalence relation).

The \emph{equivalence class} of $(p,\,q)$ is denoted by 
	\[\boxed{pq^{-1}:H\rDashTo G}\,.
\]

The situation is depicted by the following diagram:
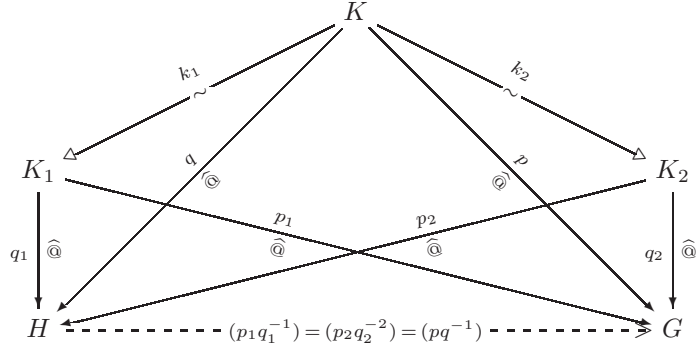
\begin{align}
	\bdi[h=3em,w=3em,tight,inline,\lab=\sst]
&&&&K&&&&\\
&&&\ldsTo^{k_1}~\sim(4,2)\ldTo^q_{\widehat{@}}(4,4)&&\rdsTo^{k_2}~\sim(4,2)\rdTo^p_{\widehat{@}}(4,4)&&&\\
K_1&&&&&&&&K_2\\
\dTo^{q_1}_{\widehat{@}}&\rdTo^{p_1\HHH}_{\widehat{@}\HHH}(8,2)&&&&&&\ldTo^{\HHH p_2}_{\HHH\widehat{@}}(8,2)&\dTo^{q_2}_{\widehat{@}}\\
H&&&&\rDashTo~{(p_1 q_1^{-1})\,=\,(p_2 q_2^{-2})\,=\,(p q^{-1})}&&&&G\\
\edi\notag
\end{align}

\begin{propdef}\label{mero}
The following properties of the pair of exactors $(p,\,q)$ are preserved by equivalence:
\begin{itemize}
	\item[(i)] $q$ is an s-equivalence;
	\item[(ii)] $p$ and $p$ are cotransversal.
\end{itemize}
When they are both satisfied, $pq^{-1}$ is called a \emph{meriedric morphism} or briefly \emph{meromorphism} from $H$ to $G$.
\end{propdef}

When these conditions are satisfied, setting $N=\text{Ker}\,p$, $R=\text{Ker}\,q$ (the latter \emph{principal}), we can write again the ``\emph{butterfly diagram}''
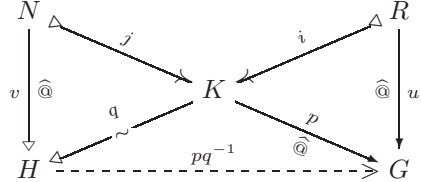
\begin{center}
\bdi[h=1.5em,w=3.5em,tight,inline,\lab=\sst]
N&&&&R\\
&\rdiTo^j(2,2)&&\ldiTo^i(2,2)&\\
\dsTo^v_{\widehat{@}}&&K&&\dTo_{u}^{\widehat{@}}\\
&\ldsTo^q~\sim&&\rdTo^p_{\widehat{@}}\\
H&&\rDashTo^{pq^{-1}}&&G
\edi
\end{center}
which owns now some added more special properties:
\begin{itemize}
	\item property (ii) is reflected by $u$ and $v$ being \emph{exactors};
	\item the (dyssymmetrical) property (i) implies $R$ being \emph{principal} and $v$ being an $s$-exactor.
\end{itemize}

From the previous section we know that $S=N\cap R$ is a smooth embedded principal subgroupoid of $K$.

\begin{propdef}\label{irred}
For a pair $(p,\,q)$ satisfying the properties of the previous proposition, the following properties are equivalent:
\begin{itemize}
	\item[(i)] $S$ is null;
	\item[(ii)] $N$ and $R$ are transverse in $K$; 
	\item[(iii)] $p$ and $q$ are cotransverse;
	\item[(iv)] $u$ is an actor;
	\item[(v)] $v$ is an actor;
	\item[(vi)] $(p,\,q)$ is a terminal object inside its equivalence class.
\end{itemize}
Then $p/q$ is called a \emph{reduced} or \emph{irreducible fraction}.
\end{propdef}

\begin{rem}
If $H$ is null ($H=E$) and $p/q$ irreducible, then $p$ is a principal actor; the orbit space of the corresponding action is the underlying space of the null groupoid $E$; $pq^{-1}$ is a \emph{non-abelian cohomology class} on $E$.
\end{rem}

Using the theory of smooth quotients of groupoids \cite{P86} to divide by $S$, we get the following basic result:

\begin{prop}\label{irredrepr}
Any meromorhism is represented by a \emph{unique irreducible fraction} with which it will be identified. In turn this irreducible representative may be identified (up to equivariant isomorphism) with a \emph{Haefliger--Skandalis--Hilsum morphism} \cite{H, HS}.
\end{prop}

The two commuting actions are defined by the $s$-actor $v$ and the principal actor~ $u$; the base of $H$ is the orbit manifold of the principal action of $G$ on the base of the groupoid $K$.

Now the use of non-irreducible representatives allows a very simple definition of the composite of two meromorphisms by means of the diagram:

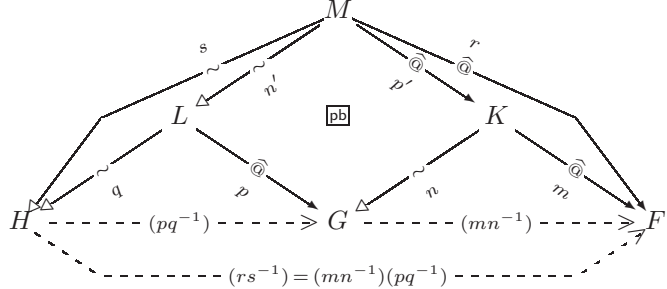
\begin{center}
\bdi[h=1em,w=3em,tight,inline,\lab=\sst]
&&&&M&&&&\\
&&&{}\ldsTo_{n'}~{\sim}(2,4)\ruLine~{\sim}^{s}(3,4)[abut]&&{}\rdTo_{p'}~{\widehat{@}}(2,4)\luLine^{r}~{\widehat{@}}(3,4)[abut]&&&\\
&&&&&&&&\\
&&&&&&&&\\
&{}&L&&\pb&&K&{}&\\
\ldsTo(1,4)[abut]&\ldsTo(2,4)~{\sim}_{q}&&\rdTo_{p}~{\widehat{@}}(2,4)&&\ldsTo~\sim_{n}(2,4)&&\rdTo_{m}~{\widehat{@}}(2,4)&\rdTo(1,4)[abut]\\
&&&&&&&&\\
&&&&&&&&\\
H&&\rDashTo~{(pq^{-1})}&&G&&\rDashTo~{(mn^{-1})}&&F\\
&\rdDashLine(1,2)[abut]&&&&&&\ruDashTo(1,2)[abut]&\\
&{}&&&\rDashLine~{(rs^{-1})\,=\,(mn^{-1})(pq^{-1})}[abut]&&&{}&\\
\edi
\end{center}
where the middle square is a pullback. By diagram chasing and a repeated use of the general properties stated in the previous sections, it can be proved that the equivalence class of the composite depends only upon the classes $(pq^{-1})$ and $(mn^{-1})$ and is again a meromorphism.

The category of meromorphisms will be denotd by $\boxed{\textbfs{\widetilde{\textsl{G}}}\textbfs{\textsl{D}}}$\,.

Now we define the (non-obvious) functor from $\textbfs{\textsl{GD}}$ to $\textbfs{\widetilde{\textsl{G}}}\textbfs{\textsl{D}}$ by means of the holograph (see above diagram (holo) in section \ref{Hologr}).

\begin{prop}\label{holograph}Let $f:H\rightarrow G$ be a smooth functor, $(p,\,q)$ its holograph.
\begin{itemize}
	\item[(i)] $p/q$ is an irreducible fraction which we \emph{identify} with the meromorphism $\widetilde{f}= pq^{-1}$.
	\item[(ii)] Two functors $f$, $g$ define the same meromorphism $\widetilde{\mathstrut f}=\widetilde{\mathstrut g}$ iff they define the same holomorphism $[f]=[g]$ \emph{(section \ref{sec:Holo})}.
	Hence we can \emph{identify} $[f]$, $\widetilde{f}$, $p/q$ and $pq^{-1}$.
	\item[(iii)] $f\mapsto\widetilde{f}$ defines a uniferous functor $\gamma=\widetilde{{\;}}:\textbfs{\textsl{GD}}\rightarrow\textbfs{\widetilde{\textsl{G}}}\textbfs{\textsl{D}}$ for which $\widetilde{\mathstrut f}\,\widetilde{\mathstrut q}=\widetilde{\mathstrut p}$, and $\gamma$ admits a factorization:
\[\boxed{\textbfs{\textsl{GD}}\stackrel{[\;]}{\longrightarrow}[\textbfs{\textsl{G}}]\textbfs{\textsl{D}}\hookrightarrow\textbfs{\widetilde{\textsl{G}}}\textbfs{\textsl{D}}}
\]
 through the canonical full functor \emph{[\;]}, by which we \emph{identify} the category of holomorphisms with a uniferous subcategory of the category of mero\-morphisms.
	\item[(iv)] A meromorphism is a \emph{holomorphism} iff it admits a representative $p/q$ with $q$
 \emph{split}. Then $v$ (in the butterfly diagram) is split too.
\end{itemize}
\end{prop}

(In particular any meromorphism with as its source a group or a plurigroup with discrete base is a holomorphism.)

\begin{thdef}\label{frac}
\begin{itemize}
	\item[(i)] The functor $\gamma:\textbfs{\textsl{GD}}\rightarrow\textbfs{\widetilde{\textsl{G}}}\textbfs{\textsl{D}}$ is the universal solution of the problem of \emph{fractions} \cite{GZ} of $\textbfs{\textsl{GD}}$ for the subcategory $\textbfs{\Sigma}$ made up with the $s$-equivalences.
	\item[(ii)] $\gamma(f)$ is an isomorphism iff $f$ is a smooth equivalence : then $\gamma(f)$
 is called a \emph{holoedric equivalence}.
 \item[(iii)] $pq^{-1}$ is an isomorphism in $\textbfs{\widetilde{\textsl{G}}}\textbfs{\textsl{D}}$ iff $p$ is an $s$-equivalence: then $pq^{-1}$ is called a \emph{meriedric equivalence}.
 \item[(iv)] The $s$-equivalences, the $i$-equivalences, the smooth equivalences and the meriedric equivalences generate the same notion of equivalence between smooth groupoids.
 \end{itemize}
\end{thdef}

The equivalence class of a smooth groupoid is therefore its isomorphy class in $\textbfs{\widetilde{\textsl{G}}}\textbfs{\textsl{D}}$. Equivalent orbital atlases are isomorphic in $\textbfs{\widetilde{\textsl{G}}}\textbfs{\textsl{D}}$.

\begin{rem}\label{GZ}
\begin{itemize}
	\item[(i)] The classical conditions for the calculus of right (nor left) fractions \cite{GZ} are \emph{not fulfilled}: we can say we have got a \emph{simplified calculus} of right fractions in that  sense that our fractions are equivalent to an irreducible (or simple) one (and are a ``right multiple'' of this one by an $s$-equivalence).
	\item[(ii)] If we identify any manifold with a \emph{null} groupoid, $\textbfs{\textsl{D}}$ \emph{is identified with a full subcategory of } $\textbfs{\widetilde{\textsl{G}}}\textbfs{\textsl{D}}$.
	\item[(iii)] The category $[\textbfs{g}]\textbfs{\textsl{D}}$ of \emph{conjugacy classes} of homomorphisms between \emph{Lie groups} is identified with a \emph{full} subcategory of $\textbfs{\widetilde{\textsl{G}}}\textbfs{\textsl{D}}$. This is valid too for plurigroups with discrete bases.
	\item[(iv)] In the case of meriedric equivalences, the butterfly diagram becomes symmetric and reversible; this special case had been presented in \cite{P77} and will be developed elsewhere: in that case the principal actors $u$ and $v$ are called \emph{conjugate}.
	\item[(v)] Given two orbital structures $Q$, $Q'$ and choosing orbital atlases $G$, $G'$ for these structures, the set $\textbfs{\widetilde{\textsl{G}}}\textbfs{\textsl{D}}(G,\,G')$ depends on the choices but up to bijections. But this does not allow to take the orbital structures for objects of a category. However this is possible when there is a \emph{canonical} choice of a meriedric equivalence between two equivalent orbital atlases: this is the case for graphoids \cite{P84} and more generally convectors in the sense of \cite{P85}.
\end{itemize}
\end{rem}

\section{Application to the fundamental group.}
\label{sec:ApplicationToTheFundamentalGroup}

In the present framework we can restate the Theorem 2 of \cite{PA} in a more striking form:

\begin{thm}
The full subcategory of discrete plurigroups is reflective \cite{ML} in th category of fractions $\textbfs{\widetilde{\textsl{G}}}\textbfs{\textsl{D}}$.
\end{thm}

In particular to any \emph{connected} orbital structure (which means the the associated topological space is connected
, there is associated a well defined (up to isomorphisms) \emph{discrete group} which, in the case of the orbital structure 
of the space of leaves of a foliation, coincides with the fundamental group of van Est-Haefliger \cite{vE, H} (and in the case of a connected smooth manifold with the Poincar\'{e} group).

This group is invariant under a wider equivalence in which uniferous retro\-connected (i.e. the fibres are connected) extensors are admitted too. This will be studied and developed elsewhere.

\providecommand{\bysame}{\leavevmode\hbox to3em{\hrulefill}\thinspace}

\appendix
\section{Comments added in March 2008}
\label{sec:CommentsAddedInJanuary2008}

\subsection{Minor modifications}
The above text is essentially a printable version (from a \LaTeX{} source) of the text \cite{P89},  not yet available in \textsc{numdam}, with the following very minor modifications (\emph{it is not at all intended to be a new up to date version}):
\begin{enumerate}
	\item The numberings of propositions may have slightly changed.
	\item The French r\'{e}sum\'{e} is translated into English.
	\item The corrections of several misprints or slips are included.
	\item Other minor more or less obvious slight limited errors are corrected.
	\item Some other minor limited modifications of the redaction are performed, aiming at a better legibility.
	\item The diagrams are redrawn, using Paul Taylor's package. The smooth equivalences, actors and ex-actors (or hyper-actors, see below) are tagged respectively by $\sim$\,, @, $\widehat{@}$, while the embeddings /surmersions are tagged by $\riTo$~/$\rsTo$, and the arrows of the category of fractions $\textbfs{\widetilde{\textsl{G}}}\textbfs{\textsl{D}}$ by $\rDashTo$.
	
	\item Several new diagrams are added (using present facilities), describing with a figure what is written in the text with words and formulas, in order to make the text easier to read (in particular the diagram (holo) of section \ref{Hologr} is presented in an equivalent but more legible way).
	\item As to the terminology:
			\begin{enumerate}
	\item ``\emph{Coarse}'' is replaced by ``\emph{banal}''\footnote{
	``Pair groupoid'' in \cite{MK'} and many other authors.}
	, which agrees with the terminology used in \cite{P07}.
	\item The use of ``\emph{Lie groupoid}'' as synonymous of ``\emph{$s$-transitive}'' is withdrawn. This agrees with the change in terminology from \cite{MK} to \cite{MK'}. Here the term ``\emph{smooth groupoid}'' is synonymous of ``\emph{differentiable groupoid}'' in the sense of Ehresmann (\cite{OC} /50/). Till 1987 the term ``Lie groupoid'' was universally used for the special case called here ``$s$-transitive''. The extension to the general case was introduced in \cite{CDW} (with the agreement of the author, consulted by P. Dazord) : it was quickly widespread from the nineties, and seems to be generally accepted nowadays. We do not use it above in order to avoid ambiguity.
			\end{enumerate}
\end{enumerate}


\subsection{Chronology.}
\label{chrono}

The existence of a connection between Haefliger--Skandalis--Hilsum generalized maps and calculus of fractions (in the sense of \cite{GZ}) in the category of groupoids, is only more or less implicit or understood in the fundamental papers \cite{H, vE, HS}. Moreover only the topological or (possibly smooth) \'{e}tale cases are considered in those papers.

The above explicit presentation was expounded orally by the author at several Conferences, notably \cite{P87, P88}, before being published in \cite{P89}, which is essentially the writing of the quoted lectures. 

Further developments and explanations about motivations were delivered in our papers quoted in the complementary bibliography below (now available on arXiv).

What was essentially new in \cite{P89} (and this will be explained in more details below) was:
\begin{enumerate}
	\item Taking into account the \emph{smooth framework} in its full generality; notably this allows, in the case of foliations a simultaneous treatment of the (canonical, non-\'{e}tale) holonomy groupoid and the (\'{e}tale as groupoids, non-canonical) transverse holonomy pseudogroups attached to various totally tranverse sub\-manifolds; but this allows also to consider smooth groupoids which are no longer (smoothly) equivalent to \'{e}tale ones, hence spaces of orbits which are much more singular than the spaces of leaves of regular foliations.
	\item The observation that the generalized H-S morphisms not only solve the general abstract G-Z problem of fractions for smooth groupoids and smooth equivalences, but they do it in a \emph{specific} and concrete way which is \emph{much simpler} than and \emph{quite different} from the (here right) standard calculus of fractions.
\end{enumerate}

\begin{enumerate}
	\item As to the first point, it should be noted that the category $\textbfs{\textsl{D}}$ of (smooth morphisms between) (possibly non-Hausdorff) smooth manifolds is very far from being a topos, and is neither complete nor cartesian closed; it demands quite different methods developed in \cite{P89}.
	
	As explained in more details for instance in \cite{P07}, the point of view which is systematically adopted (in view of wide-ranging generalizations) consists first in an ``\emph{internalization}'' of various algebraic or set-theoretic definitions and constructions by means of suitable \emph{diagrams}, in which certain pullbacks, certain monomorphisms (here possibly non-proper embeddings), and certain epimorphisms (here the surmersions) are stressed, and second, rather than looking for actual limits (which in general won't exist), considering certain suitable \emph{equivalence classes} of such diagrams (which indeed encapsulate a richer information than actual limits).
	
	For instance, the following algebraic notions are in this way carried over to the smooth framework, in which they immediately make sense:
\begin{itemize}
	\item injectivity (for a mapping);
	\item fullness and faithfulness (for a functor $f$), described by means of the diagram $\textbf{T}(f)$ (see diagram (T) page \pageref{AT});
	\item surjectivity (for a mapping);
	\item essential surjectivity (for a functor) (Def. \ref{essur});
	\item equivalence, here between groupoids (see Prop.Def. \ref{ind})\footnote{
	This definition and the previous one appear in \cite{M02} without any reference. This notion is sometimes (notably in \cite{M88} in the étale case, and in \cite{Mr}) named ``essential equivalence'', while the term of equivalence is kept for the internalization of adjoint equivalence \cite{ML}, a very unfortunate and unacceptable terminology, since the former notion is weaker (and also more fundamental and less special) than the latter. In order to repare this incoherency, several authors introduced later the term ``weak equivalence'', which is equally unfortunate, since, though weaker than an adjoint (differentiable) equivalence, this notion is stronger than a mere equivalence (in the purely algebraic sense): one has to \emph{add} conditions of differentiability required by the internalization.
	};%
	\item action laws, discrete (op)fibrations, (op)fibrations (in the categorical, not topological, sense) (see diagram (A) page \pageref{AT}, and Def. \ref{act});
	\item principal actions and bundles; H-S bibundles (see Prop.Def. \ref{mero}, and \ref{irred} and Prop. \ref{irredrepr}; the set-theoretic description of these bi\-bundles is encapsulated in the ``butterfly diagram'').
\end{itemize}

\item As to the second point (to be developed below in full details) it is important to notice that:
\begin{itemize}
	\item the precise conditions for G-Z calculus of (right) fractions are \emph{not} fulfilled (even in the \'{e}tale or topological cases), neither in the category $\textbfs{\textsl{G}}\textbfs{\textsl{D}}$, nor in the category $\textbfs{[\textsl{G}}]\textbfs{\textsl{D}}$ (see Dictionary below (\ref{dico}))\footnote{
	Contrary to what is claimed in \cite{M02} and \cite{L01}, as explained in full details below.
	};%
		\item anyway this question is immaterial since the general description of fractions of \cite{GZ} cannot yield the much simpler and precise description one gets here;
	\item \emph{not all pairs} $(p,\,q)$ (with $q$ a smooth equivalence) are considered (see Prop.Def. \ref{mero}) but only those for which $p$ is an ex-actor (see Dictionary below (\ref{dico}), and which furthermore satisfy a certain \emph{co-trans\-versality condition}; in particular this condition \emph{cannot be satisfied when $q$ is a unit} (except if $p$ is a unit too).
	\item on the other hand the equivalence between such pairs $(p,\,q)$ which is considered here is \emph{stricter} than the G-Z equivalence;
	\item very remarquably such an equivalence class admits a canonical \emph{irreducible representative} of the fraction (see Prop. \ref{irredrepr}), which may be precisely \emph{identified with a H-S generalized morphism}, and the other representatives are just \emph{right multiples} (by an $s$-equivalence) of the irreducible one, exactly as for integers and rational numbers (save for the lack of commutativity), whence our terminology of ``simplified calculus of fractions''; note that, though these other representatives don't admit a simple set-theoretic bi-bundle interpretation by means of two commuting action laws (one of them being principal), however their \emph{diagrammatic description} (using ``actors'' embedded in a ``butterfly diagram'') is nearly as simple, just replacing actors by ex\-actors or hyper\-actors (see Dictionary (\ref{dico}));
	\item the same simple formal construction describes at the same time, as categories of fractions, the categories $\textbfs{[\textsl{G}}]\textbfs{\textsl{D}}$ and $\textbfs{\widetilde{\textsl{G}}}\textbfs{\textsl{D}}$, just considering in one case the split $s$-equivalences, and in the other the surmersive equivalences, and it turns out, very remarkably, that all the smooth equivalences are then inverted (see Prop. \ref{holograph} and Th.Def. \ref{frac})\footnote{
	This contrasts with the method suggested in \cite{M02} and \cite{L01}, which uses the category $\textbfs{[\textsl{G}}]\textbfs{\textsl{D}}$ as an intermediate step towards $\textbfs{\widetilde{\textsl{G}}}\textbfs{\textsl{D}}$, in which the G-Z method would apply; even if this method might work, it would be anyway a roundabout and very awkward way for describing fractions.
	}; %
	but the use of sur\-mersive equivalences instead of starting directly with general smooth equivalences leads to over\-simplifications\footnote{
	We remind that, even from a purely algebraic poin of view, the equi\-valences are \emph{not stable by pulling back}; that seems to create insuperable difficulties when trying to check the G-Z conditions.
	};
	\item one gets also the fact that $\textbfs{[\textsl{G}}]\textbfs{\textsl{D}}$ is faithfully embedded in $\textbfs{\widetilde{\textsl{G}}}\textbfs{\textsl{D}}$\footnote{
	Very strangely this seems in \cite{L01} to be considered implicitely as a general property (used in an essential way in the sequel of the paper) of the G-Z calculus of fractions, which is of course obviously wrong.
	};%
	\item though this is not fully exploited in \cite{P89}, but just suggested, the same purely diagrammatic method applies to various interesting sub\-categories of the category of all smooth equivalences (for instance considering proper surmersions, or surmersions with connected or simply connected fibres, see \cite{P07}), and also larger categories, such as our sub\-category of ``extensors'', and various intermediate ones; this deserves further study and may have many applications.
\end{itemize}
\end{enumerate}

One finds also in \cite{P89} precise criteria for the existence of fibred products; a special paragraph will be devoted to that point below.

\subsubsection{}

It has to be noted that the other authors listed below do not quote \cite{P89}, save \cite{Mr}, where it is quoted but not used; in this latter paper the three frameworks of \'{e}tale, topological and general smooth groupoids and morphisms are so inextricably mixed that it is often hard to understand why and whether some definitions or statements are expressed in one of these frameworks rather than another.

However we intend to show below in detail that:
\begin{itemize}
	\item some of these papers rediscover, in a more or less precise way, a part of the assertions of \cite{P89};
	\item in all of them the precise link between the G-Z fractions and the H-S generalized morphisms is \emph{totally misunderstood}.
\end{itemize}

\subsection{Dictionary}
\label{dico}
In order to make comparisons easier we mention here some correspondences in terminology between the papers quoted in the complementary bibliography below.
\subsubsection{Later changes in terminology}
\label{newterm}

	In \cite{P07} some \emph{changes in terminology} are introduced:
\begin{itemize}
	\item concerning the (commutative) \emph{squares}, the terms :
\begin{center}
\emph{universal}, \emph{$i$-faithful}, \emph{$s$-faithful}
\end{center}
are replaced by:
\begin{center}
\emph{gpb, ipb, spb}
\end{center}
(where gpb means ``good pullback'');
\item as to the \emph{morphisms}, the terms:
\begin{center}
\emph{ex}/\emph{in}-\emph{actor}.
\end{center}
are replaced by: 
\begin{center}
\emph{hyper}/\emph{hypo}-\emph{actor}
\end{center}
\end{itemize}
We remind that in the purely algebraic setting (i.e. forgetting the smooth structures), the current terminology, in the categorical literature, for our ``\emph{hyper}-\emph{actors} /\emph{actors}'' is 
\begin{center}
``(\emph{op}\-)\emph{fibrations} /\emph{discrete} (\emph{op})\-\emph{fibrations}'',
\end{center}
 but these terminologies cannot be extended to the topological and smooth contexts, in which they have quite different well-established meanings, which creates an unacceptable ambiguity. 

[Note that the same is true about the use, for groupoids, of 
\begin{center}
``\emph{discrete}\-/\emph{coarse}\-/\emph{connected}'',
\end{center}
which is widespread in the categorical literature, where we use instead:
\begin{center}
``\emph{null}\-/\emph{banal}\-/\emph{transitive}''.)]
\end{center}
Moreover it is clear that the most significant part of the content of these terms (in the categorical context) degenerates and becomes trivial when leaving large categories for small groupoids, for which just a certain surjectivity condition remains.

\subsubsection{Other authors.}
\begin{itemize}
	\item \emph{Categories of fractions}: the correspondence for the names of the categories of fractions in \cite{P89,L01,M02} is given by

\begin{align}
&\textnormal{\cite{P89}}&&\textnormal{\cite{L01}}&&\textnormal{\cite{M02}} \notag\\
	&[\textbfs{\textsl{G}}]\textbfs{\textsl{D}}&&\textsf{LG'}&&\textsl{G}_0 \notag\\
	&\textbfs{\textsl{\widetilde{G}}}\textbfs{\textsl{D}}&&\textsf{LG}&&\textsl{G}_0[W^{-1}] \notag
\end{align}

\item \emph{Smooth equivalences}: see footnote $(^2)$ above in \ref{chrono}.

\end{itemize}

We now enter into more details in order to justify our claims.

		\subsection{Calculus of right fractions revisited.}
			\label{GZcalcul}

In order to be in a position to stress the specific aspects of the ``\emph{simplified calculus of fractions}'' of \cite{P89} (briefly called HSH calculus, since it uses generalized morphisms in the sense of Haefliger-Skandalis-Hilsum), we start by reminding some basic facts about the classical calculus of (here right) fractions as expounded in \cite{GZ} (called GZ calculus), adopting a slightly more general viewpoint. 

For the convenience of the reader the \emph{left} calculus of \cite{GZ} will be translated (\emph{by dualization}) into the language of \emph{right} calculus involved here.

\subsubsection{``Concrete'' versus ``abstract'' fractions.}
\label{abscon}

First, for an arbitrary category $\textbfs{\textsl{C}}$, the problem of formally inverting the arrows of a given subcategory (or even subclass) $\textbfs{\Sigma}$ (such arrows are tagged here by $\rsTo^{\simeq\h}$, by analogy with our $s$-equivalences $\rsTo^{\sim\h}$) has 
always\footnote{It has been noticed in later literature that some logical caution is necessary with respect to the ``sizes'' of the universes involved, but that is not the point here.}
 a universal ``abstract non-sense'' solution $\textbfs{\textsl{\widetilde{C}}}=\textbfs{\textsl{C}}[\textbfs{\Sigma}^{-1}]$. 
 
However the arrows of this latter category are (somewhat unformally) described by classes of chains of arrows (to be read from left to right) of the form 
\begin{align}
	\cdot\lsTo^{\simeq}_{q_1} \cdot\rTo_{p_1} \cdot \h\cdots\h\cdot\lsTo^{\simeq}_{q_n} \cdot\rTo_{p_n} \cdot
	\notag
\end{align}
 up to an obvious equivalence relation: the sub-chains of type 
\begin{center}
 \bdi[w=2em,tight,abut,\lab=\sst,inline]
\,H\,&\lsTo_{q}^{\simeq}&\,K\, &\rsTo_{q}^{\simeq} &H\\
\edi 
\h or 
 \bdi[w=2em,tight,abut,\lab=\sst,inline]
\,H\,&\rsTo_{q}^{\simeq}&\,K\, &\lsTo_{q}^{\simeq} &H\\
\edi 
\end{center}
have to be cancelled.

But such classes are not easy to handle concretely. This will be referred to in the sequel as the ``\emph{abstract}'' solution (of the \emph{problem of fractions} for $\textbfs{\Sigma}$) and an equivalence class of diagrams of the previous form will be referred to below as an \emph{abstract fraction}.

Now what will be called below a ``\emph{concrete}'' solution consists in looking for a \emph{description} of 
the new arrows by \emph{much simpler diagrams}, which are just certain pairs of the form 
\begin{align}
(p,\,q)=
\left(
\bdi[w=2em,tight,abut,\lab=\sst]
\,H\,&\lsTo_{q}^{\simeq}&\,K\, & \rTo~{p} &G\\
\edi
\right)\,,\notag
\end{align} 
\emph{up to a certain suitable equivalence relation}. Such an equivalence class will be referred to as a \emph{concrete \emph{(right)} fraction}.

\emph{It is important to notice that, in such a general presentation, we do not demand all the pairs $(p,\,q)$ to be considered, nor the concrete fraction to coincide (as an equivalence class) with the whole abstract fraction} (though this will be the case for the GZ calculus, but not for our simplified HSH calculus).

It results that, though the abstract solution is essentialy unique, \emph{there may exist different concrete solutions} in the sense just defined.

Now one finds in \cite{GZ} \emph{sufficient conditions}, and the description of the suitable \emph{equivalence relation} between pairs $(p,\,q)$, for the existence of such a \emph{concrete solution}.

\subsubsection{ G-Z conditions.}
\label{GZcond}
\begin{conv}\label{conv}
\emph{From now on it is understood that in the figures below the dotted arrows and the brown symbols will denote arrows or objects which are to be constructed in the text, or the existence of which is asserted in the stated axiom or condition, contrasting with the initial data.}
\end{conv}

The G-Z conditions for the calculus of \emph{right} fractions are \emph{dual} to the conditions listed (a), (b), (c), (d) in \cite{GZ} for the left fractions, and may be denoted by adding a star ${}^{\ast}$.

Conditions (a) and (b), which are self-dual, just mean that $\textbfs{\Sigma}$ is a uniferous\footnote{
i.e. having the same objects or units. This is sometimes called in the literature large or ample subcategory, which is not very suggestive. The terminology we use here agrees with the one used for rings with unit.
} %
sub\-category of $\textbfs{\textsl{C}}$, and there is no problem about them.

Conditions (c) and (d) are stated in a minimal form in \cite{GZ}. 

\begin{itemize}
	\item Condition $(\textnormal{c}^\ast$) is implied 
by the stronger condition (which may be sometimes easier to check):

	$(\textnormal{C}^\ast$) every pair \begin{align}
(f,\,s)=
\left(
\bdi[w=2em,tight,abut,\lab=\sst,inline]
\,H\,&\rsTo_{s}^{\simeq}&\,G\, & \lTo~{f} &G'\\
\edi
\right)\,,\notag
\end{align} 
can be completed into a pullback square of the following type\footnote{This is the property named ``parallel transfer'' in \cite{P07}.}:

\begin{center}
\bdi[h=1.8em,w=4em,tight,inline,textflow]
\emphor{H'\pbk}&\rTo^{\emphor{f'}}[dotted]&H\\%
&&\\
\dTo^{\emphor{s'}}[dotted]~{\emphor{\simeq}}&&\dsTo^{s}~{\simeq}\\%
&&\\
G'&\rTo^f&G\\
\edi
\end{center}
(condition $(\textnormal{c}^\ast$) demands just the square to be commutative).

This strong condition is satisfied in $\textbfs{\textsl{G}}\textbfs{\textsl{D}}$ by our $s$-equivalences, but \emph{not} when taking directly \emph{all} the (smooth) eqivalences. We shall examine below the attempt of working in $[\textbfs{\textsl{G}}]\textbfs{\textsl{D}}$ proposed in \cite{M02}.

\item condition $(\textnormal{d}^\ast$) is expressed briefly by the following commutative diagram (in which the data are in black, while the existence of the brown object and the dotted arrows  is asserted by the condition):

\bdi[h=3em,w=3em,p=0.5em,tight]
\emphor{D}&\rsTo~{\emphor{\simeq}}[dotted]&A&&\\
&\rdTo(2,2)[dotted]&\dTo\dTo&\rdTo(2,2)&\\
&&B&\rsTo~{\simeq}&C\\
\edi

One may observe that this condition is trivially verified when all the arrows of $\textbfs{\Sigma}$ are monomorphisms, and cannot be verified (save in very degenerate situations) when they all are epimorphisms.

This condition cannot be verified when taking for $\Sigma$ our $s$-equivalences.%

\end{itemize}

\subsubsection{GZ versus HSH fractions }
\label{sec:GZVersusHHHFractions}

\begin{enumerate}
	\item In GZ calculus one considers \emph{all} pairs $(p,\,q)$ with denominator $q\in\Sigma$, while in section \ref{sec:Mero}, one demands moreover the numerator $p$ to be an ``ex-actor'', and the pair to be ``co-transversal''.
	\item The equivalence between such pairs is not the same for GZ and HSH calculus and may be described picturially by figure \ref{fig:defeq} below (for HSH the arrows $k_i$ have to belong to $\Sigma$).

\begin{figure}[htbp]
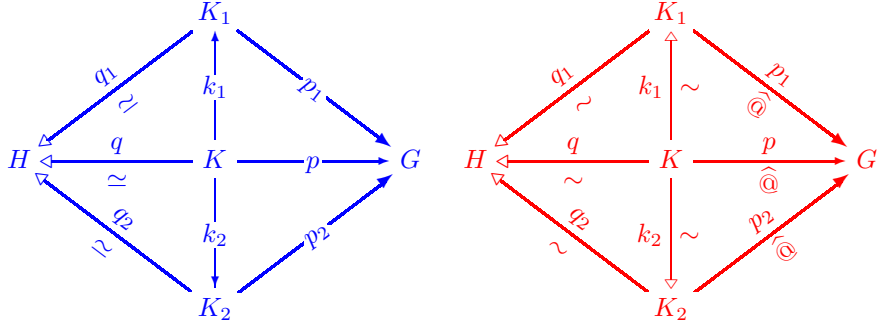

\begin{center}
\emphb{\bdi[h=2.8em,w=3.7em,tight,inline]
&&K_1&&\\
&\ldsTo^{q_1}_\simeq[thick]&\uTo~{k_1}&\rdTo~{p_1}[thick]&\\
H&\lsTo^q_{\simeq}&K&\rTo~p&G\\
&\lusTo^{q_2}_\simeq[thick]&\dTo~{k_2}&\ruTo~{p_2}[thick]&\\
&&K_2&&\\
\edi}\H
\emphr{\bdi[h=2.8em,w=3.7em,tight,inline]
&&K_1&&\\
&\ldsTo^{q_1}_\sim[thick]&\usTo^{k_1}_{\sim}&\rdTo^{p_1}_{\widehat{@}}[thick]&\\
H&\lsTo^q_{\sim}&K&\rTo^{p}_{\widehat{@}}&G\\
&\lusTo^{q_2}_\sim[thick]&\dsTo^{k_2}_{\sim}&\ruTo^{p_2}_{\widehat{@}}[thick]&\\
&&K_2&&\\
\edi}
\end{center}
	\caption{Definition of the equivalence between two fractions $(p_i,\,q_i),\;(i=1,\,2)\,$: 
	\emphb{GZ calculus of right fractions} versus \emphr{HSH} \emphr{simplified calculus}.}
	\label{fig:defeq}
\end{figure}

The proof of the transitivity of this equivalence, using GZ conditions, though not so much obvious, is left to the reader in \cite{GZ}. In figure \ref{fig:trans} below (where Conv. \ref{conv} is used) it is suggested by the blue diagram (one gets first $M$, then $L''$), while the red one suggests the simpler proof when condition (C$^\ast$) is satisfied (one gets  $L''$ directly by taking the fibred product of the arrows $k_2$ and $k'_2$); however, for HSH calculus, one has here to prove also the co-transversality of ($p_3,\,q_3$), which needs more care.

An analogous simplification occurs in the definition of the \emph{composition} of the fractions (section \ref{sec:Mero}).

\begin{figure}[htbp]
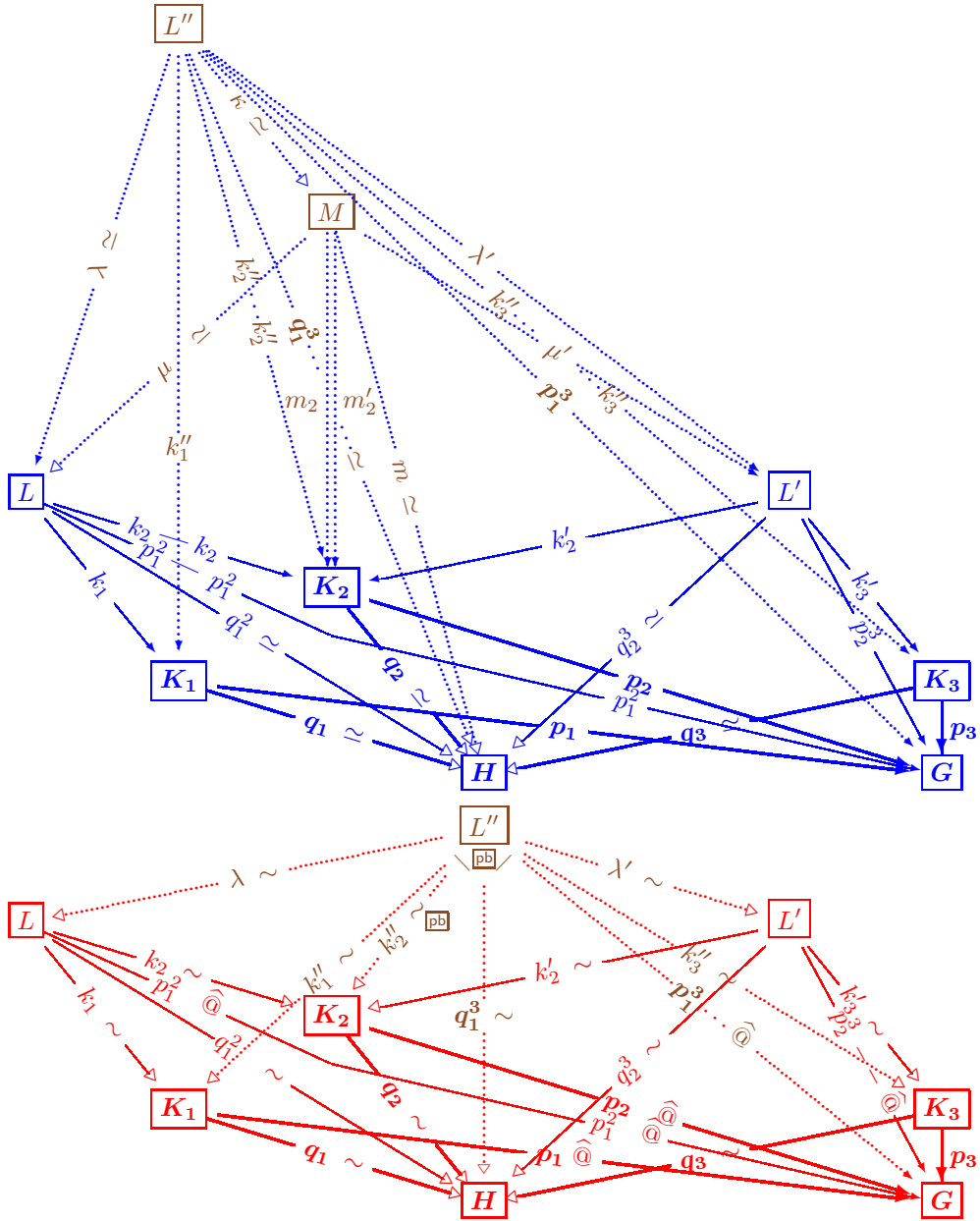

\emphb{\bdi[h=1.8em,w=2.94em,p=0.3em,tight,\lab=\textstyle]
&&\emphor{\boxed{L''}}&&&&&&&&&&\\
&%
\ldTo~{\emphor{\;\lambda\;\simeq\;}}(2,10)[dotted]
&&%
\rdsTo~{\emphor{\;\textbfs{q_1^3}\emphb{\cdots\H\cdots}\emphor{\simeq}\;}}(4,16)[dotted]
\rdsTo~{\emphor{\;\kappa\;\simeq\;}}(2,4)[dotted]
\rdTo~{\emphor{\;\lambda'\;}}(8,10)[dotted]
\rdTo~{\emphor{\;\textbfs{p_1^3}\;}}(10,16)[dotted]
\rdTo~{\emphor{\;k''_2\;\emphb{\cdots} k''_2\,}}(2,12)[dotted]
\rdTo~{\emphor{\;k_{3}''\emphb{\cdots\h\h\cdots} k_{3}''}}(10,14)[dotted]
&&&&&&&&&\\
&&&&&&&&&&&&\\
&&&&&&&&&&&&\\
&&&&\emphor{\boxed{M}}&&&&&&&&\\
&&&
\ldsTo~{\emphor{\;\mu\h\simeq\;}}(4,6)[dotted]
&&
\rdTo~{\emphor{\;\mu'\;}}(6,6)[dotted]
\rdsTo~{\emphor{\;m\;\simeq\;}}(2,12)[dotted][abut]
&&&&&&&\\
&&&&&&&&&&&&\\
&&&&&&&&&&&&\\
&&&&
\dTo^{\emphor{m_2}}[dotted][abut]
\dTo_{\emphor{m_{2}'}}[dotted][abut]
&&&&&&&&\\
&&\dTo~{\emphor{k_{1}''}}[dotted]&&&&&&&&&&\\
\boxed{L}&&&&&&&&&&\boxed{L'}&&\\
&%
\rdsTo~{\;q_1^2\;\simeq\;}(6,6)
\rdTo~{\;k_1\;}(2,4)
\rdTo~{\;k_2\;\emphb{\textbf{---}}\; k_2\;}(4,2)
\rdLine~{\;p_1^2\;\emphb{\textbf{---}}\;\; p_1^2\;}(4,3)[abut]
&&&&&&&&\ldTo~{\;k'_2\;}(6,2)
\ldsTo~{\;q_2^3\;\simeq\;}(4,6)
&&\rdTo~{\;k'_3\;}(2,4)
\rdTo~{\;p_2^3\;}(2,6)
&\\
&&&&\boxed{\textbfs{K_2}}&&&&&
&&
&\\
&&&&{}&
\rdsTo~{\;\textbfs{q_2}\;\simeq\;}(2,4)[thick][abut]
\rdTo~{\;\textbfs{p_2}\;}(8,4)[thick]
&&&&&&&\\
&&\boxed{\textbfs{K_1}}&&&\rdTo~{\;p_1^2\;}(8,3)[abut]&&&&&&&\boxed{\textbfs{K_3}}\\
&&&%
\rdsTo~{\;\textbfs{q_1}\;\simeq\;}(4,2)[thick][abut]
\rdTo~{\;\textbfs{p_1}\;}(10,2)[thick]
&&&&&&&&\ldsTo~{\;\textbfs{q_3}\;\simeq\;}(6,2)[thick][abut]
&%
\dTo_{\textbfs{p_3}}[thick][abut]\\
&&&%
&&&\boxed{\textbfs{H}}&&&&&
&\boxed{\textbfs{G}}\\
\edi}
\emphr{\bdi[h=1.8em,w=2.94em,p=0.3em,tight,\lab=\textstyle]
&&&&&&\emphor{\un{\diagdown\pb\diagup}{\boxed{L''}}}&&&&&&\\
&&&&&%
\ldsTo~{\emphor{\;\lambda\;\sim\;}}(6,2)[dotted]
\ldsTo~{\emphor{\,k''_2\;\sim\,}}(2,4)[dotted]
\ldsTo~{\emphor{\;k_{1}''\;\sim\;}}(4,6)[dotted]
&&%
\rdsTo~{\emphor{\;\lambda'\;\sim\;}}(4,2)[dotted]
\rdsTo~{\emphor{\;k_{3}''\;\sim\;}}(6,6)[dotted]
\rdTo~{\emphor{\textbfs{p_1^3}}\;\emphr{\cdots}\;\emphor{\widehat{@}}\;}(6,8)[dotted]&&&&&\\
\boxed{L}&&&&&{\H\;\emphor{\pb}}&&&&&\boxed{L'}&&\\
&%
\rdsTo~{q_1^2\h\sim\;}(6,6)
\rdsTo~{\;k_1\;\sim\;}(2,4)
\rdsTo~{\;k_2\;\sim\;}(4,2)
\rdLine~{\,p_1^2\h\widehat{@}\,}(4,3)[abut]
&&&&&
&&&\ldsTo~{\;k'_2\;\sim\;}(6,2)
\ldsTo~{\;q_2^3\;\sim\;}(4,6)
&&\rdsTo~{\,k'_3\;\sim\;}(2,4)
\rdTo~{\,p_2^3\;\---\;\,\widehat{@}}(2,6)
&\\
&&&&\boxed{\textbfs{K_2}}&&\dsTo~{\emphor{\;\textbfs{q_1^3}\,\sim\;}}[dotted]&&&
&&
&\\
&&&&{}&
\rdsTo~{\;\textbfs{q_2}\;\sim\;\;}(2,4)[thick][abut]
\rdTo~{\,\textbfs{p_2}\h\widehat{@}\,}(8,4)[thick]
&&&&&&&\\
&&\boxed{\textbfs{K_1}}&&&\rdTo~{\,p_1^2\h\widehat{@}\,}(8,3)[abut]&&&&&&&
\boxed{\textbfs{K_3}}\\
&&&%
\rdsTo~{\;\textbfs{q_1}\;\sim\;}(4,2)[thick][abut]
\rdTo~{\;\textbfs{p_1}\;\widehat{@}\;}(10,2)[thick]
&&&&&&&&\ldsTo~{\;\textbfs{q_3}\;\sim\;}(6,2)[thick][abut]
&%
\dTo_{\textbfs{p_3}}[thick][abut]\\
&&&%
&&&\boxed{\textbfs{H}}&&&&&
&\boxed{\textbfs{G}}\\
\edi}
	\caption{Proof of transitivity for the equivalence between two fractions: 
	\emphb{GZ calculus of right fractions} versus \emphr{HSH simplified calculus}.}
	\label{fig:trans}
\end{figure}

\item The basic simplification in the description of HSH fractions comes from Prop. \ref{irredrepr}, which states that they are just right multiples of an irreducible representative (identified with an HSH bi-bundle) by an s-equivalence.

This basic result cannot be obtained by the GZ calculus.

\end{enumerate}

\subsection{HSH calculus in $\protect\textbfs{\textsl{G}}\textbfs{\textsl{D}}$ versus  GZ calculus in $\protect[\textbfs{\textsl{G}}]\textbfs{\textsl{D}}$}

Another aspect of the oversimplification obtained in section \ref{sec:Mero} as compared to the strategy sketched in \cite{M02,L01} and others is that one works directly in $\textbfs{\textsl{G}}\textbfs{\textsl{D}}$ and one starts taking as $\textbfs{\Sigma}$ the subcategory of $s$-equivalences, which has very nice properties of stability, especially by pulling back, and are much easier to handle than the general equivalences, and then Prop. \ref{holograph} yields an alternative description of 
$[\textbfs{\textsl{G}}]\textbfs{\textsl{D}}$; \emph{a posteriori} the subcategory of general smooth equivalences appears as the saturation (in the sense of \cite{GZ}) of the subcategory of $s$-equivalences.

\subsubsection{Working in $\protect[\textbfs{\textsl{G}}]\textbfs{\textsl{D}}$}
The method of \cite{M02} and others would need a preliminary study of $[\textbfs{\textsl{G}}]\textbfs{\textsl{D}}$ in order to be able to check whether the conditions of GZ calculus (not fulfilled in $\textbfs{\textsl{G}}\textbfs{\textsl{D}})$ are now satisfied. (We note that on the contrary it would be much easier to study this category using its alternative description as a category of fractions, which inverts \emph{split} $s$-equivalences, as given by Prop. \ref{holograph} and Th.-Def. \ref{frac}).

Such a study is made conscientiously in \cite{Pk}, using 2-categories and precise conditions of coherence, demanding much work and care, but, though this is achieved in an apparently very general and abstract framework, this can be applied only to the étale case and not to the general differentiable framework.

This latter framework is approached in \cite{M02} in a very sketchy way, where it is proposed to use diagrams which are commutative up to homotopy, which amounts to forget or neglect the coherence conditions; moreover the same notation is used for the arrows in $\textbfs{\textsl{G}}\textbfs{\textsl{D}}$ and their images in $[\textbfs{\textsl{G}}]\textbfs{\textsl{D}}$, which creates an inextricable confusion. As was done above we shall avoid this confusion by using below dashed arrows for the arrows of $[\textbfs{\textsl{G}}]\textbfs{\textsl{D}})$.

In this prospect \cite{M02} introduces a fibred product in $[\textbfs{\textsl{G}}]\textbfs{\textsl{D}}$ ambiguously called ``fibred product'' and denoted in the same way as the fibred product in $\textbfs{\textsl{G}}\textbfs{\textsl{D}}$.

To avoid this ambiguity, we call it here \emph{weak} fibred product and denote it by $\widetilde{H'}=G'\un{G}{\widetilde{\times}}H$.

The precise relation between the (true) fibred product and the weak one is given by the diagram of Fig. \ref{fig:pb} below, where one can recognize the occurrence of the holographs introduced in section \ref{Hologr} (cf. diagram (holo)).

This relation is obtained by pulling back along the canonical $i$-equivalence $\iota_g:G\rightarrow \sq G$, and this yields an $i$-equivalence $\epsilon$ between the two fibred products. Note that $\epsilon$ admits an inverse in $\widetilde{\textbfs{\textsl{G}}}\textbfs{\textsl{D}}$, but not in $[\textbfs{\textsl{G}}]\textbfs{\textsl{D}}$.


One can use this diagram and our Prop. \ref{pb} to study the properties of existence and stability of this weak product in $[\textbfs{\textsl{G}]\textbfs{\textsl{D}}}$. One sees that, as claimed in \cite{M02}, the general equivalences become now stable by pullback \emph{when this pullback does exist}. But this existence still requires a certain condition of transversality which is not always satisfied (in spite of the essential surmersivity).

For this reason condition (C$^\ast$) above is not satisfied in $[\textbfs{\textsl{G}]\textbfs{\textsl{D}}}$ when taking all the equivalences.

For the same reason it does not seem possible to prove \emph{directly} the transivity of the Morita equivalence as defined in \cite{M02}, i.e. using directly for the \emph{definition} all the equivalences\footnote{
Anyway it cannot be proved just using the weak fibred product, as claimed in \cite{M02}.
}, though this definition is \emph{a posteriori} correct as a \emph{consequence} of one of the main \emph{theorems} of the theory developed in \cite{P89}, where one starts with the $s$-equivalences, easy to handle (cf. Th.-Def. \ref{frac} (iv)).

Anyway, from the remarks above, \emph{this question is strictly immaterial and uninteresting} since, even if these conditions are satisfied, \emph{the GZ construction don't yield here any useful information about the description of the arrows of the category} 
$\widetilde{\textbfs{\textsl{G}}}\textbfs{\textsl{D}}$, which is already known to \emph{exist} by the general abstract nonsense argument.

This will be our final conclusion.

\begin{figure}[htbp]
\emph{\begin{diagram}[h=1.5em,w=2.14em,tight]
\emphr{\boxed{\emphr{\textbfs{H'}}}}{\negthickspace\emphr{\pbk}}
&&\rTo~{\emphr{\boxed{\emphr{h}}}}[thick]
&&&&\emphr{\boxed{\emphr{\textbfs{H}}}}&&&&&&&&&&\\
&\rdiTo~{\emphvi{\boxed{\emphvi{\varepsilon}}\,\sim\,}}(4,4)[thick][abut]
&&&&
\ruDashTo~{\Bigl[{\overline{h}}\Bigr]}(2,4)[thick]
&&\rdiTo~{\,\eta\,\sim\,}(4,4)[abut] 
 \rdTo~{\emphvi{=\boxed{\emphvi{\text{Id}_H}}=}}(10,4)[thick]
 &&&&&&&&&\\
\dTo~{\emphr{\boxed{\emphr{u'}}}}[thick]&&&&&&
\dLine~{\emphr{\boxed{\emphr{u}}}}[thick]
&&&&&&&&&&\\
&&&&&{}&\HonV&&\hLine~{\emphb{\boxed{\overline{h}}}}[abut][thick]&{}&&%
\hLine~{\emphb{\widehat{@}}}[abut][thick]%
&&\HonV&&{}&\\
&&&&\ruLine(1,1)[abut][thick]
\emphb{\boxed{\emphb{\textbfs{\widetilde{H'}}}}}{\pbk}
&%
&\HonV&&
\rTo~{h^{\square}}
&&H^{\square}\pbk
&&&\rsTo~{q(u)\,\sim\,}
&&&\emphb{\boxed{\emphb{\textbfs{H}}}}
\rdTo[abut][thick](1,1)\\
&&&\ruLine(1,1)[abut][thick] 
\ldDashTo~{\Bigl[\overline{u'}\Bigr]}(4,3)[thick]
&&\rdDashTo(2,3)[abut]
\rdTo(6,14)[abut]
\rdTo(12,7)[abut]
&&&&\ruLine(1,1)[abut]
&&&&&&\ldDashTo~{\h\;}^{\HHHHH\H [u]=p(u)/q(u)}(6,14)[thick]&\\
&&&&\dTo~{u'^{\square}}&&\dTo[thick]
&&&&\dTo~{\;u^{\square}}&&&&&&\dTo~{\emphb{\boxed{u}}}[thick]\\
\emphr{\boxed{\emphr{\textbfs{G'}}}}
&&\rLine~{\emphr{\boxed{\emphr{g}}}}[thick]
&\VonH&\VonH&
\rTo[thick]&
\emphr{\boxed{\emphr{\textbfs{G}}}}&&&
&&&&&&&\\
&\rdiTo~{\,\gamma\,\sim\,}(4,4)[abut] \rdTo~{\emphvi{=\boxed{\text{Id}_{G'}}=}}(4,11)[thick]&&&&&& \rdiTo~{\emphv{\boxed{\iota_G}\,\sim\,}}(4,4)[thick][abut] \rdTo~{\emphvi{=\boxed{\text{Id}_G}=}}(4,11)[thick] \rdTo~{\,\emphvi{=\boxed{\text{Id}_G}=}\,}(10,4)[thick]
&&&
&&&&&&\\
&&&\dLine~{\ov{{}}{\emphb{\boxed{\overline{u'}}}}}[abut][thick]&&&&&&\dLine~{p(u)}[abut]&&&&&&&\\%
&&&&&{}&&&\hLine~{{p(g)}}[abut]&{}&&
\hLine~{\widehat{@}}[abut]%
&&&&{}&\\
&&&{}
&\ruLine(1,1)[abut]G'^{\square}\pbk
&&&&\rTo~{g^{\square}}&&\emphv{\pmb{\sq G}}
&&&\rsTo~{\emphv{\boxed{\varpi_2}\,\sim\,}}[thick]
&&&\emphb{\boxed{\emphb{\textbfs{G}}}}\rdTo(1,1)[abut]\\
&&&&&&&&&%
%
&&&&&&
\ruDashTo(6,7)~{\emphb{={\boxed{\emphb{\Bigl[\text{Id}_G\Bigr]}}}=}}[thick] \ldDashTo(6,7)~{{\phantom{=\boxed{\emphb{\Bigl[\text{Id}_G\Bigr]}}=}}}[thick] 
\ruDashTo~{\h}^{[g]=p(g)/q(g)\HHHHH\HH}(12,7)[thick]&\\
&&&\dLine~{\ov{{}}{\emphb{\widehat{@}}}}[abut][thick]&&&&&&\dLine~{\ov{}{\widehat{@}}}[abut]
&&&&&&&\\
&&&&{}\dsTo~{\h q(g)\sim\,}
&{}&&&&&
\dsTo~{\emphv{\boxed{\varpi_1}\,\sim\,}}[thick]&&&&&&\\
&&&&&&&&&&&&&&&&\\
&&&&&&&&&&&&&&&&\\
&&&{}&&&&&&{}&&&&&&&\\
&&&&\rdTo(1,1)[abut][thick]
\emphb{\boxed{\emphb{\textbfs{G'}}}}
&&&&\rTo~{\emphb{\boxed{g}}}[thick]&&
\rdTo(1,1)[abut]
\emphb{\boxed{\emphb{\textbfs{G}}}}&&&&&&\\
\end{diagram}}
The diagram of full arrows is (strictly) commutative in $\textbfs{\textsl{GD}}$.\\
When adding the dashed arrows, it is still (strictly) commutative in $[\textbfs{\textsl{G}}]\textbfs{\textsl{D}}$. \\
All the squares are (good) pullbacks (in $\textbfs{\textsl{GD}}$).
	\caption{Strict pullback \emphr{$H'=G'\un{G}{\times}H$} vs weak pullback \emphb{$\widetilde{H'}=G'\un{G}{\widetilde{\times}}H$}.}
	\label{fig:pb}
\end{figure}

\providecommand{\bysame}{\leavevmode\hbox to3em{\hrulefill}\thinspace}

\end{document}